\documentclass[onecolumn,aps,floatfix,superscriptaddress]{revtex4}
\usepackage{amsmath,amssymb,eucal,graphicx}
\usepackage{amsbsy}
\usepackage{graphicx}
\usepackage{epsfig}
\usepackage{comment}
\usepackage{array}
\usepackage{xcolor}
\usepackage[colorlinks=true, urlcolor=blue, anchorcolor=blue, citecolor=blue,filecolor=blue,linkcolor=blue]{hyperref}

\newcommand{\Real}{\mathbb{R}}
\newcommand{\Natural}{\mathbb{N}}
\newcommand{\bfx}{\boldsymbol{x}}
\newcommand{\bfp}{\boldsymbol{p}}

\newcolumntype{P}[1]{>{\centering\arraybackslash}p{#1}}

\newtheorem{thm}{Theorem}

\newtheorem{corollary}{Corollary}
\newtheorem{assume}{Assumption}
\newtheorem{rem}{Remark}

\definecolor{myred}{rgb}{0.8, 0.1, 0.1}

\begin{document}
\title{Knowledge-informed neuro-integrators for aggregation kinetics}
\author{D. Lukashevich}
\affiliation{Skolkovo Institute of Science and Technology, Moscow, 121205, Russia,\\ e-mail: dmitrii.lukashevich@skoltech.ru}
\author{I. Tyukin}
\affiliation{King's College London, London, Strand, WC2R 2LS, United Kingdom, \\ e-mail: ivan.tyukin@kcl.ac.uk}
\author{N. Brilliantov}
\affiliation{Skolkovo Institute of Science and Technology, Moscow, 121205, Russia}
\affiliation{School of Computing and Mathemarical Sciences, University of Leicester, Leicester LE1 7RH, UK, \\ email: nb144@leicester.ac.uk}

\begin{abstract}
We report a novel approach for the efficient computation of solutions of a broad class of large-scale systems of non-linear ordinary differential equations, describing aggregation kinetics. The method is based on a new take on the dimensionality reduction for this class of equations which can be naturally implemented by a cascade of small feed-forward artificial neural networks.
We show that this cascade, of otherwise static models,  is capable of predicting solutions of the original large-scale system over large intervals of time,  using the information about the solution computed over much smaller intervals. The computational cost of the method depends very mildly on the temporal horizon,   which is a major improvement over the current state-of-the-art methods, whose complexity increases super-linearly with the system's size and proportionally to the simulation time. In cases when prior information about the values of solutions over a relatively small interval of time is already available, the method's computational complexity does not depend explicitly on the system's size. The successful application of the new method is illustrated for spatially-homogeneous systems, with a source of monomers, for a number of the most representative reaction rates kernels. 


\end{abstract}

\maketitle

\section{Introduction}

Aggregation phenomena are ubiquitous in nature and utilized in numerous industrial processes \cite{Krapivsky,Leyvraz}. These processes occur at many spatial and temporal scales, ranging from molecular to astrophysical ones. Aggregation at molecular scales may be exemplified by classical processes, such as irreversible polymerization and agglomeration of colloidal particles \cite{Davis2004}. Aggregation of prions (biological macromolecules) is an example of these processes occurring in living cells; it is associated with multiple  Alzheimer-like diseases \cite{prions}. Other prominent examples in the realm  of biological systems include coagulation of erythrocytes (single blood cells) to form roleaux \cite{EP1926,RWS1982,RWS1984},  aggregation of bacteria via dextran induction \cite{VR1980}, formation colonies of viruses \cite{CoagVirus} and, on much larger scales, -- of schools of fish \cite{CoagFish}. Outside of biology, aggregation occurs in the atmospheric processes, when airborne particulates  coalesce into smog particles  \cite{Srivastava1982}, aggregation of water droplets underlying clouds formation \cite{Falkovich2002}, etc. At the astrophysical scale, agglomeration of particles is an important process, driving the formation of planetary rings and planetesimals \cite{PNAS,Cuzzi,Esposito_book}. Aggregation processes are also widespread in networks of different nature, including social, economic \cite{CoagNetw}, and internet communities \cite{Krapivsky,Dorogov}. In the context of networks, aggregation and fragmentation correspond to the nucleation, merging, and splitting of different forums of users \cite{Dorogov}. 

Importantly,  modelling aggregation and fragmentation processes involves the necessity to compute solutions of large-scale or even infinite-dimensional systems of coupled ordinary differential equations known as the Smoluchowski equations \cite{Krapivsky,Leyvraz}. The computational complexity of state-of-the-art methods, such as the low-rank decomposition methods \cite{matveev2015, matveev2018parallel,skorych2019investigation,osinsky2020low,MKSTB,OscPRE2017}, scales as $O((N\log N) N_T)$, where $N$ is the system's size and $N_T$ is the temporal horizon of the modelling expressed by the number of temporal steps. The challenge, however, is that we are often interested in understanding the behaviour of systems with large $N$ over long intervals of time. This requires significant computational resources.  Given the practical importance of this class of equations, there is hence a need for the development of more computationally efficient methods and tools.

In this work, we show that combining ideas from model reduction with tools developed in the seemingly distant area of neural networks enables the derivation of approximate solutions of the Smoluchowski equations at the computational costs scaling as $O((N\log N) N_\tau) + O(m N_T)$, $N_\tau \ll N_T$, $m\ll N$. This is a significant improvement over the best conventional solvers when the modelling time horizons $N_T$ are large. The workflow of the proposed solver is as follows: 

\begin{itemize}

\item At the first step, we compute the solution of the original initial value problem over a short interval of time $[0,\tau]$, $\tau>0$, using any traditional solver. Invoking available knowledge about general properties of the solutions (e.g. non-negativity, monotonic decay, etc.), we select an appropriate set of functions, parameterized by a time-dependent vector of parameters, capable of approximating solutions of the Smoluchowski equations across all relevant scales at any fixed time in the interval of interest $[0,T]$, $T>0$, $T\gg\tau$. 

\item At the second step, we employ an artificial neural network (ANN) to implement the parameterized dependence established at the first step, with parameters now being the weights and biases of the network. We call this ANN a "parameterizing NN". For any fixed value of time, the parameterizing NN, subject to the choice of its parameters, is capable of approximating solutions of the Smoluchowski system. However, to be able to use this network for computing solutions of the original system at all relevant values of time, one needs to have an additional process determining the values of weights and biases of the parameterizing ANN which would correspond to these values of time. Since we are interested in computing solutions of the Smoluchowski system beyond the original short interval of time, this process must be able to extrapolate knowledge contained in the values of the solutions we have already obtained using the traditional solver at the first step of the method.

\item To address this problem, at the third step, we identify appropriate properties of the time evolution of the parameters vector of the parameterizing NN and formulate the requirements imposed on the extrapolation. Then we construct the second ANN; its function is to predict the values of parameters of the parameterizing NN over long times. We call this ANN "predictive NN". 

\item Finally, in the fourth step, we apply our predictive NN to produce an approximation of the solution of the initial value problem for the original system of the aggregation equations over the entire interval of time of interest. 

\end{itemize}

One of the remarkable and unexpected outcomes of our work is that it shows that the dynamics of large-scale, and potentially infinite-dimensional, Smoluchowski systems with kernels for which analytical solutions are not yet known can be described by mere ten nonlinear first-order differential equations. Moreover, they can be implemented by a cascade of two very small ANNs: (i) a parameterizing network with two neurons and (ii) a predictive network with twelve neurons resulting in the relatively low computational costs of the method. We would also like to note that the method can potentially be applied to a fully data-driven setup in which the initial step of computing the solution over the interval $[0,\tau]$ is replaced with acquiring or measuring empirical data over that interval. In this case, the computational complexity of the method reduces to $O(m N_T)$ (see Section \ref{sec:discussion} for more detials).

The rest of the paper is organized as follows. In  the next Section \ref{sec:formal} we present a brief technical introduction into modelling of aggregation kinetics, introduce the Smoluchowski equations, provide a statement of the problem and present a formal framework describing the proposed method. In Section \ref{sec:prediction} we present the application of the new method to a class of equations with the most representative aggregation kernels. In Section \ref{sec:discussion} we discuss the limitations of the method, produce estimates of its computational costs, and outline the potential directions of its future development. Section \ref{sec:conclusion} summarizes our findings and concludes the paper.

\section{Knowledge informed dimensionality reduction for aggregation processes}\label{sec:formal}

 In this section, we present a formal description of the proposed approach to efficiently compute solutions of Smoluchowski equations. We begin with describing all relevant features of the process and briefly introduce the Smoluchowski equations. This is then followed by the introduction of the formal approach to the problem and a description of  the method, involving neural networks.

\subsection{Aggregation processes}

\begin{figure}[h]
\begin{center}
\includegraphics[scale=0.21]{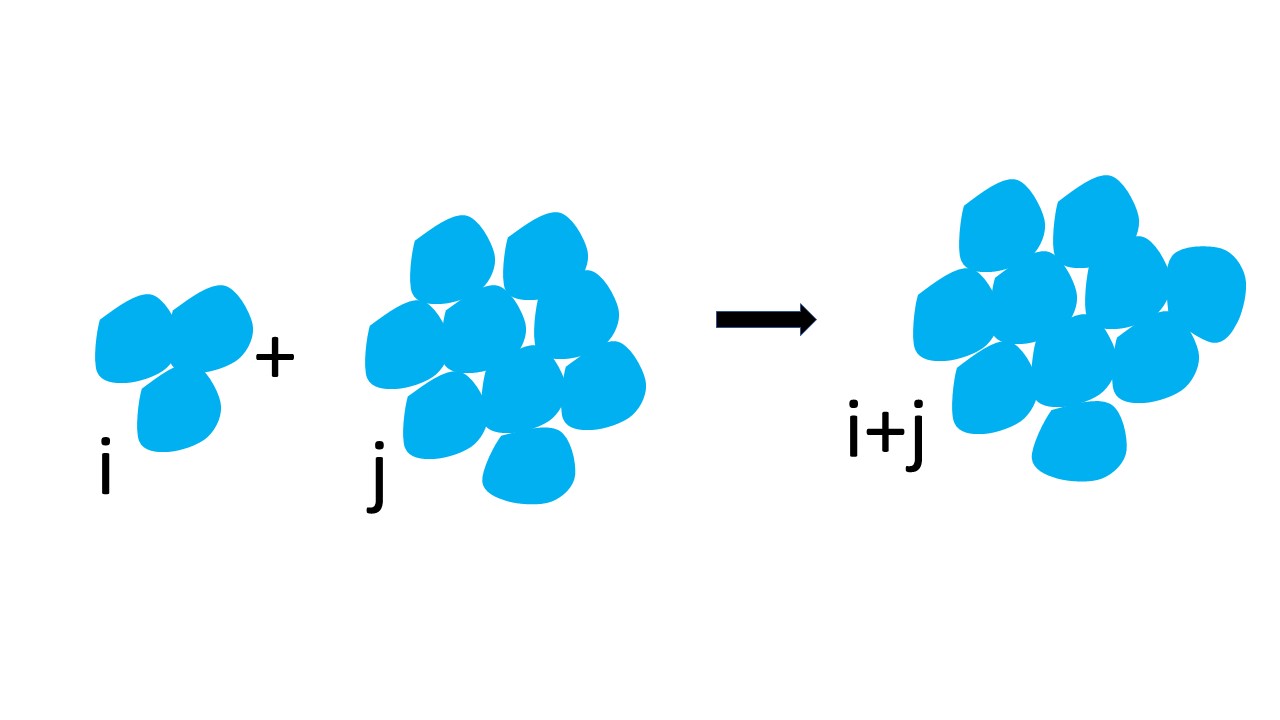}
\end{center}
\caption{An aggregate of size $i$ (that is, comprising $i$ monomers) meets with an  aggregate of size $j$ to form a joint aggregate of the size $i+j$. }
\label{fig:Agg_Frag}
\end{figure}

The aggregation takes place when two objects, comprising respectively $i$ and $j$ elementary objects (monomers), merge and form a joint aggregate of $i+j$ monomers (see Fig.~\ref{fig:Agg_Frag}). Formally, it can be expressed by the kinetic equation
$$
[i] +[j] \xrightarrow{K_{ij}} [i+j],
$$
where $K_{ij}$ is the merging rate. 

The above process can be generalised to systems of arbitrarily large size. For such large systems, a common and natural practice is to consider densities (or concentrations) instead of the sets of individual aggregating objects. Let $c_k$ be the densities of aggregates of size $k$, that is, of aggregates comprising $k$ monomers. The rate of change of $c_k$ is determined by Smoluchowski equations \cite{Krapivsky,Leyvraz}

\begin{equation}
 \label{eq:Smol1} \frac{dc_k}{dt} = \frac{1}{2}\sum_{i+j=k}K_{i,j}c_i c_j - c_k
\sum_{i=1}^{\infty}K_{i,k}c_i, ,\qquad k = 1..\infty   
\end{equation}
The first term in the right-hand side of the above equation describes the formation rate of $k$-mers from clusters of size $i$ and $j$, while the second term accounts for the loss of $k$-mers due to agglomeration of these clusters with all other aggregates. The factor $1/2$ in the first term is used to exclude  double counting of the same process ($i+j \to k$ and $ j+i \to k$). The important property of the kinetic equations \eqref{eq:Smol1} is the conservation of the total mass of all clusters: 
\begin{equation}
\label{eq:Mass} M(t)=\sum_{k=1}^{\infty} k c_k(t) = {\rm const.}
\end{equation}
Equations \eqref{eq:Smol1} describe spatially-homogeneous systems. The kernels $K_{i,j}$ may be obtained from the microscopic analysis of the aggregation processes, see e.g. \cite{Krapivsky,PNAS,BrilBodKrap2009,BFP2018}. In applications, $K_{i,j}$ are usually homogeneous functions of the masses $i$ and $j$ of merging clusters. Here we will focus on the most widely used kernels of this type \cite{Krapivsky,Leyvraz}: 

\begin{equation}
\label{eq:Cijgen}
 K_{i,j} = i^{\nu}j^{\mu} + i^{\mu}j^{\nu}, \ \mu, \nu\in\Real.  
\end{equation} 

Note that Smoluchowski equations \eqref{eq:Smol1} form an {\it infinitely large } set of coupled {\it non-linear} ODE of the first order, whose solution,  ${\bf C}(t)  \equiv  (c_1(t), c_2(t), \dots)$, is an  infinitely large vector of densities determined for all $t>0$. Generally, Eqs. \eqref{eq:Smol1}-\eqref{eq:Cijgen} are analytically intractable. Only a small class of kernels admits solutions whose analytic expressions are known. This includes kernels with $K_{ij}=A+B(i+j)+Cij$, where $A,\, B, \, C$ are some constants \cite{Krapivsky,Leyvraz}, and the addition kernels \cite{bk}. In general,  one needs to use numerical methods to solve these equations. Obviously, only a finite set of equations may be treated numerically. Therefore, the following natural question arises: to what extent one can describe the infinite set of equations with a finite set and how large this finite set must be?

For conventional solvers, the minimal number of equations depends on the initial conditions and on the intended temporal modelling horizon. Initial conditions associated with practical applications commonly define the number of non-zero densities of small aggregates, $c_1(0), \, c_2(0), \, \dots c_l(0)$ with $l$ less than $5$ whilst assuming that the density of all other aggregates are zero. Among these, the most natural and important initial condition is the mono-disperse initial condition, $c_1(0)=1$ and $c_k(0)=0$ for $k \geq 2$ (in convenient units for the density).  As time goes by, larger and larger aggregates appear requiring more and more equations to describe their evolution. If $k_{\rm max}(T)$ is the maximal size of clusters at time $T$, then the number of equations $N(T)$ which  accurately describe the evolution of the system over the interval $[0,T]$ should satisfy the condition: $N[T] > k_{\rm max}(T)$. This is the consequence of the conservation of mass property \eqref{eq:Mass} which must hold true at $t=T$. Scaling analysis, see e.g. \cite{Leyvraz}, shows that $k_{\rm max}(T)$ rapidly grows with $T$, often, as a power law, $\sim T^{a}$, with $a>1$. As a result, a very large number of equations is needed to describe the evolution of the system even over a relatively short interval of time.  This creates a vicious cycle: with the increasing number of equations, $N (T)$,  the computational complexity of deriving numerical solutions of the aggregation equations, using conventional solvers, grows rapidly, and the problem quickly becomes intractable.  

To overcome the problem, several numerical approaches have been developed. The most successful among these are methods based on low-rank decompositions of the rate kernels. These methods approximate $K_{ij}$ by a combination of kernels of small rank $R$, yielding a significant reduction of the computational complexity up 
to $\sim (R N \log N )\, N_T$, instead of $\sim (N^2) \, N_T $, 
inherent for the straightforward integration approach, see e.g. \cite{matveev2015, matveev2018parallel,skorych2019investigation,osinsky2020low,MKSTB,OscPRE2017}. Still, as the amount of computations grows rapidly with $T$, the problem remains which limits the application of conventional methods for large $T$. Moreover, not all kernels allow accurate low-rank decompositions.

As we show below, a potential approach to address this problem can be found at the interface of the ideas stemming from the area of dimensionality reduction (e.g. \cite{DimRed}) with the rapidly developing field of artificial neural networks  (e.g. \cite{AIbook}).  The paradigm of dimensionality reduction (DM) suggests that the description of a large system could be achieved, up to some acceptable accuracy, by tracing the evolution of a much smaller, reduced system. Such a reduced system may be a set of variables, parametrizing an appropriate class of functions. The choice of the class of functions typically depends on  available phenomenological knowledge about the evolution of the system. This motivates the name of the approach:  "knowledge-informed". As we show below, the exploitation of ANN allowseffective optimization of the corresponding set of parameters, yielding good approximations of the solution in the chosen class of functions. 

In practice, this approach allows predicting the evolution of large systems of non-linear aggregation equations \eqref{eq:Smol1} over large intervals of time by solving a much smaller system and, as we show here, over much shorter time intervals. 

\subsection{Truncated aggregation equations}

When an appropriate solution of the original system \eqref{eq:Smol1} is lacking, the remaining option is to proceed with solving a truncated, finite-dimensional system of $N$ equations, 
\begin{equation}
    \label{eq:Smol-trunc}
    \frac{d \hat{c}_k}{dt} = \frac{1}{2} \sum_{i+j=k} K_{ij} {\hat{c}}_i \hat{c}_j - \hat{c}_k \sum_{j=1}^{N} K_{kj} \hat{c}_j,\qquad k = 1,\dots,N, 
\end{equation}
where $\hat{c}_1(t),\ldots,\hat{c}_N(t)$  are the approximations of the densities $c_1(t),\dots,c_N(t)$, corresponding to the infinite system  \eqref{eq:Smol1}. 
The value of $N$ there depends on the duration of the evolution time $T$ needed to be explored. Certainly, the choice of $N$ in the truncated system affects the accuracy of the approximation of $c_l(t)$ by $\hat{c}_l(t)$, for $l=1, \dots, N$. Indeed, removing higher masses from consideration changes the second term in the right-hand-side of (\ref{eq:Smol-trunc}) at every $t$ and leads to the violation of the mass conservation. The impact of this change grows with time resulting in the progressing loss of accuracy for large $t$.

Even more challenging, with respect to the requirement of mass conservation, are the aggregation equations with a source -- usually, a source of monomers. Such equations are very important for numerous practical (especially industrial see e.g. \cite{MBE,Sire}) applications. In this case Smoluchowski equations read \cite{Krapivsky,Leyvraz}: 
\begin{align}
    \label{problem:smol-bin}
    & \frac{d \hat{c}_k}{dt} = \frac{1}{2} \sum_{i+j=k} K_{ij} \hat{c}_i \hat{c}_j - \hat{c}_k \sum_{j=1}^{N} K_{kj} \hat{c}_j + \delta_{k,1},\qquad k = 1,\dots,N, 
\end{align}
with the source of monomers of unit intensity, in the physical units chosen here. We assume the most important for applications mono-disperse initial conditions, $\hat{c}_k(0) = \delta_{k,1}$, where $\delta_{i,j}$ is the Kronecker delta, which immediately implies $M(t)=t+1$. In what follows we will illustrate the application of our new method to the system with the monomer source, Eqs. \eqref{problem:smol-bin}. 

\subsection{General framework}

We begin with a high-level description of our approach for computing solutions of large-scale aggregation equations.  As we will show later, the approach enables overcoming the barrier of the exploding computational complexity of the aggregation equations with $N$.  Technical details, appropriate assumptions, justification of the feasibility of the method, and error bounds expressed in terms of these assumptions are presented in the Appendix.

At the heart of the approach is the idea to exploit the knowledge of the values of ${\bf \hat{C}}(t) =\left(\hat{c}_1(t),\hat{c}_2(t),\dots, \hat{c}_{N}(t) \right)$ over some appropriately small interval $[0,\tau]$, $\tau\ll T$ for the reconstruction of
${\bf C} (t)=\left( c_1(t),c_2(t),\dots, c_{N}(t), \ldots \right)$ for all $t$ in much larger interval $[0,T]$, $T>0$. 
We expect that our method would be particularly transparent and illustrative for aggregation equations with homogeneous kernels \eqref{eq:Cijgen}, as prior theoretical studies \cite{Krapivsky,Leyvraz} indicate the existence of some scaling functions capturing the structure of solutions of the Smoluchowski equations with such kernels. 
The method assumes that the distribution of densities  $c_1(t),\dots,c_k(t),\dots$, as a function of size $k$ at a given time instance $t$, can be modelled by a pair of parameterized functions (scaling functions)

\begin{equation}\label{eq:family_of_distributions}
F: \mathbb{N}\times \mathbb{R}^m \rightarrow \mathbb{R}_{\geq 0}, \ {\bf p}: \mathbb{R}_{\geq 0} \rightarrow \mathbb{R}^m
\end{equation}
such that the following holds true
\begin{equation}\label{eq:matching}
F(k,{\bf p}(t))=\hat{c}_k(t) + \varepsilon_k(t),
\end{equation}
where 
\[
\|\varepsilon_k(t)\| \leq \Delta \ \mbox{for all} \ k\in \mathbb{N}, \ k\leq N, \ t\in [0,T] \subseteq \mathbb{R}_{\geq 0}
\]
is an approximation error whose norm is bounded by $\Delta\in\mathbb{R}_{\geq 0}$. Moreover, we will further request that the evolution of ``parameter'' vector ${\bf p}(t)$ obeys a finite system of ordinary differential equations:
\begin{equation}\label{eq:parameters_system}
\dot{\bfx}=f({\bfx},\theta), \ {\bf p}(t)=h({\bfx}(t;\theta, {\bfx}_0),\theta), \ {\bf x}_0\in\mathbb{R}^n,
\end{equation}
where $f:\mathbb{R}^n\times\mathbb{R}^d\rightarrow \mathbb{R}^n$ is a piece-wise continuous and locally Lipschitz function, $h:\mathbb{R}^n\times\Real^d \rightarrow \mathbb{R}^m$ is a continuous locally Lipschitz function, $\theta\in\mathbb{R}^d$ is a vector of parameters,  $\bfx_0\in\mathbb{R}^n$ is the vector of initial conditions, and ${\bfx}(\cdot;\theta,{\bfx}_0)$ is the maximal solution of the initial value problem (\ref{eq:parameters_system}). We assume that the domain of this solution's definition includes the target interval $[0,T]$. 

Equations (\ref{eq:family_of_distributions}) -- (\ref{eq:parameters_system}) do not merely aim to capture the behavior of the distribution ${\bf \hat{C}}(t)$ as a function of $t$. Our main goal is to {\it extrapolate} ``knowledge'' of ${\bf \hat{C}}(t)$ for $t\in[0,\tau]$ to a much larger interval of interest $[0,T]$ with $T\gg \tau$. That is why we embed a system of ODEs, represented by (\ref{eq:parameters_system}), into the core of the problem definition. Thanks to the existence and uniqueness of solutions of (\ref{eq:parameters_system}), the setting defined by (\ref{eq:family_of_distributions})--(\ref{eq:parameters_system}) 
justifies the feasibility of the task to simultaneously achieve a finite-dimensional reduction of the original system (\ref{problem:smol-bin}) and produce a computationally efficient {\it extrapolation} of estimates $\hat{c}_k(t)$ beyond the interval $[0,\tau]$ circumventing the need to solve the original computationally-demanding problem (\ref{problem:smol-bin}) over $[\tau,T]$.

Notwithstanding the potential advantages of the above approach for computing and predicting behavior of the original infinite-dimensional model, the problem of determining parameters and initial conditions of (\ref{eq:parameters_system}) may still be computationally challenging. Here we propose to overcome this issue by combining phenomenological insights with some peculiar properties of a class of neural networks with sigmoidal activation functions. 

\subsection{The class of functions and parametrization of the solution. Parametrizing network}

\begin{figure}[ht!]
    \centering
    \includegraphics[width=0.9\columnwidth]{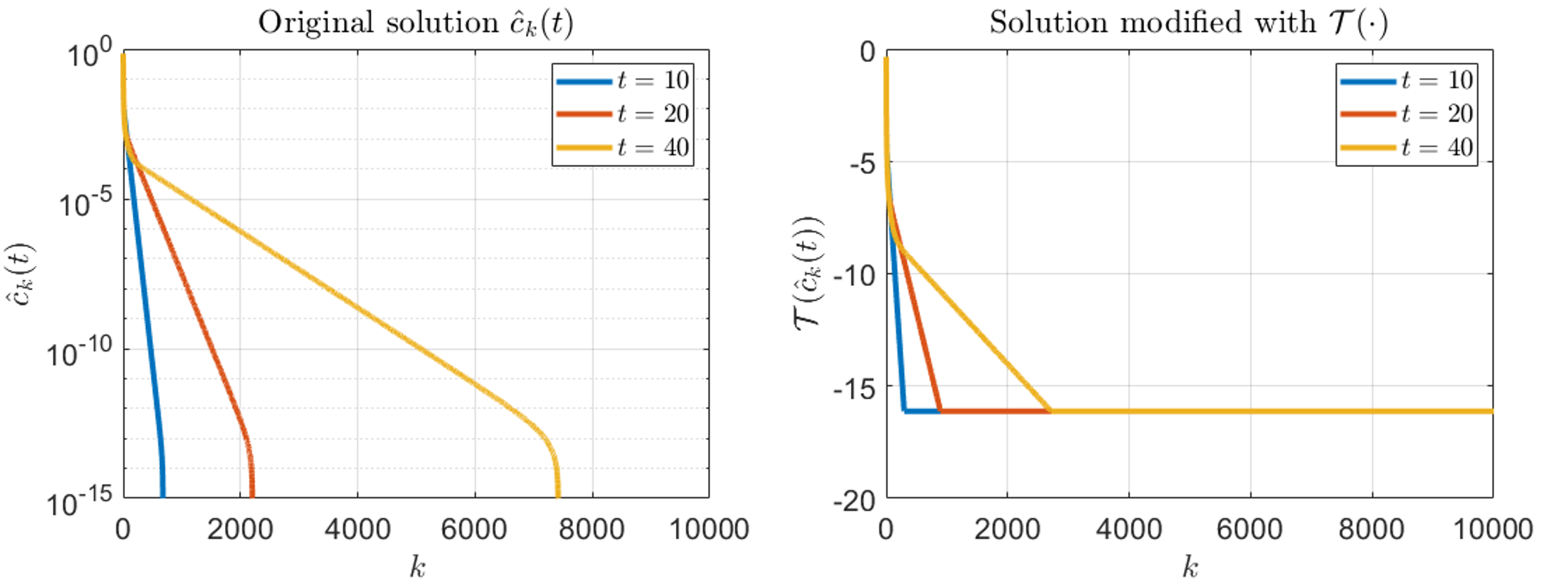}
    \caption{\label{fig:orig-trans-solution} Original solution of the truncated system (\ref{problem:smol-bin}) with (\ref{monodisp}) depicted with logy scale (left) and modified with transformation (\ref{eq:solution-transformation}) (right).}
\end{figure}

Fig. \ref{fig:orig-trans-solution} (left panel) depicts the solution of the truncated system (\ref{problem:smol-bin}). As certain values of the solution are below machine precision, they cannot be properly simulated and thus utilized. To avoid gaps in data we assume that such values are equal to a small constant $\varepsilon = 10^{-7}$. 

The other relevant point to discuss is the potential advantage of introducing an endogenous scaling inside the map $F$. The motivation behind such scaling is to balance the values of the distribution $\hat{c}_k$ as a function of $k$ between two end-points corresponding to $k=1$ and $k=N$. The values of $c_k$ can range from $0$ to $1$, and as the size of the system grows, many variables $c_k$ may  take values close to $0$. When this happens, any algorithm's fitting and predicting capabilities could be hindered by the low sensitivity in this region. That is why, in this paper, we propose to apply a logarithmic endogenous scaling, implemented by the mapping $\mathcal{T}$, as defined by (\ref{eq:solution-transformation}) below. The effect of the application of the mapping $\mathcal{T}$ to $\hat{c}_j$, where
\begin{align}
    \label{monodisp}  
   & \hat{c}_k(0) = \delta_{k,1} \\ 
    \label{Kij1}
   & K_{ij} = 1 
\end{align}
is illustrated with Fig. \ref{fig:orig-trans-solution} (right panel). In what follows we define $\mathcal{T}$ as:
\begin{equation}
    \label{eq:solution-transformation}
    \mathcal{T} (\hat{c}_k(t)) = \begin{cases}
        \log (\hat{c}_k(t)), \mbox{  }  \hat{c}_k(t) > \varepsilon\\
        \log (\varepsilon), \mbox{  } \hat{c}_k(t) \leq \varepsilon.
        \end{cases}
\end{equation}

If we find a mapping $\mathcal{F}$ such that
\[
\mathcal{F}(k,{\bf p} (t))=\mathcal{T}(\hat{c}_k(t)) + \epsilon_k(t),
\]
then the original parameterised mapping $F$ introduced in (\ref{eq:family_of_distributions}) can be defined as
\[
F(k,{\bf p}(t)) = \left\{\begin{array}{ll}\exp \left(\mathcal{F}(k,{\bf p}(t))\right), &  \mathcal{F}(k,{\bf p} (t)) > \log(\varepsilon)\\
            \varepsilon, &  \mathcal{F}(k, {\bf p}(t)) = \log(\varepsilon)
            \end{array}\right.
\]
with $\varepsilon_k(t)$ in (\ref{eq:matching}) accounting for the propagation of the error term through the transformation. Once the mapping $\mathcal{F}(k, {\bf p}(t))$ is found, the problem of extrapolation of the solution $\mathcal{T}(\hat{c}_k(t))$ can be reduced to the problem of extrapolation of the parameter vector ${\bf p}(t)$.

The mappings $\mathcal{F}$ can be taken from a broad class functions, including splines, polynomial, or piece-wise linear functions. Here we will look for $\mathcal{F}$ in the class of functions $NN_p$ defined as follows:
\begin{equation}
    \label{eq:NN_p_equation_gen}
    NN_p = L^m \sigma L^{m-1} \sigma \dots \sigma L^1,
\end{equation}
where $\sigma$ is a coordinate-wise activation function, and $L^i$ is an affine map
\[
    L^ix = W^i x + b^i,
\]
where $W^{i}\in\Real^{D_{i} \times D_{i-1}}$, $b^{i}\in\Real^{D_i}$ are the matrices of weights and biases, and the set $\mathbf{D} = (D_m, D_{m-1}, \dots, D_1, D_0) \in \mathbb{N}^{m+1}$ is a set of dimensions.

The class $NN_p$ defines a class of feedforward neural networks with $m$ layers \cite{higham2019deep}. Our choice of this class is motivated by two factors. First, neural networks are well-known for their universal approximation capabilities \cite{devore2021neural}. Second, there is a wide range of numerical tools and open-access codes which could be used to fit models from this class to data.

The parametrization vector ${\bf p}$ in general case reads, ${\bf p} =\left(W^m, W^{m-1}, \dots, W^1, b^m, b^{m-1}, \dots, b^1 \right)$, where all instances in the expression are vectorized. Given the relative simplicity of the curves in Fig. \ref{fig:orig-trans-solution} one can significantly simplify the above networks. In this work, we will start with the following network family,
\begin{equation}
    \label{eq:NN_p_equation}
    NN_p(k; t) = W(t) \cdot ReLU \left[ wk + B(t) \right] + b,
\end{equation}
where the values of scalars $W(t)$ and $B(t)$ are learnable, and $w = -1$ and $b = \log(10^\varepsilon)$ are frozen, $t$ is a hyperparameter representing the  time point $t$  at which the network $NN_p(k;t)$ approximates $\hat{c}_k(t)$. $ReLU(\cdot)$ is the rectified linear unit, an activation function defined as
\[
ReLU(x) = \left\{\begin{array}{ll}x, & x > 0 \\
            0, &  x \leq 0.
            \end{array}\right.
\]

\begin{figure}[ht!]
    \centering
    \includegraphics[width=0.7\columnwidth]{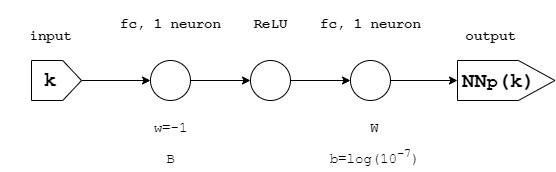}
    \caption{\label{fig:NN_p_structure} Structure of the parametrizing neural network.  It generates a parametric family of functions that contains the modified solution $\mathcal{T}(\hat{c}_k(t))$}
\end{figure}
\begin{figure}[ht!]
    \centering
    \includegraphics[width=0.7\columnwidth]{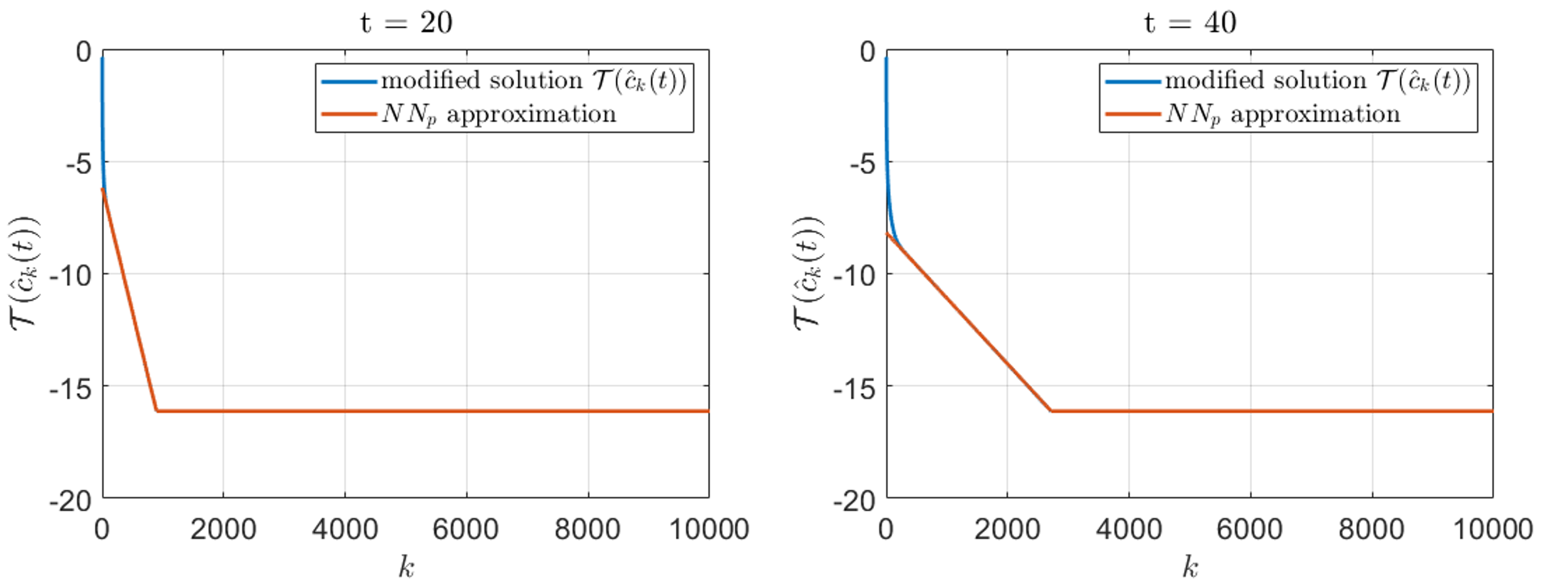}
    \caption{\label{fig:NN_p_approx} Approximation of the modified solution of the system (\ref{problem:smol-bin}) with (\ref{monodisp}) by a ReLU network $NN_p$ on different time steps.}
\end{figure}

In this context, the network itself serves as the family $\mathcal{F}(k,p(t))$, while its parameters $W(t)$ and $B(t)$, when considered as functions of time, constitute the vector-function ${\bf p}(t)=\left(W(t),B(t) \right)$ in (\ref{eq:parameters_system}).  We call this network ”parametrizing network”. As one can see from Fig. \ref{fig:NN_p_approx}, this very simple network $NN_p$ has the potential to accurately reproduce the solution $\mathcal{T}(\hat{c}_k(t))$. 

The fact that a simple ReLU-network can approximate the solution of the Smoluchowski system (\ref{problem:smol-bin}) is not much of a surprise. Available analytical solutions for such systems are exponentially decaying functions of $k$ in the scaling regime, e.g. \cite{Leyvraz}. Such a dependence, coupled with our choice to cut off all values below machine precision, is in an excellent agreement with the expressivity of this small ReLU network and its capability to approximate  $\mathcal{T}(\hat{c}_k(t))$, see Fig. \ref{fig:NN_p_approx}. What is surprising however,  is that the above parametrization, as we shall show later, enables to predict the evolution of the potentially infinite-dimensional system (\ref{problem:smol-bin}) using the evolution of a two-dimensional vector of parameters implemented by another static network. 

To predict the values of $\hat{c}_k(t)$  in $ (\tau, T]$, one needs to extrapolate the values of  $W(t)$ and $B(t)$, available over the interval $ [0,\tau]$,  to the much larger time interval $(\tau,T]$.  The respective extrapolation should account for the most important properties of these functions; eventually, this will result in an accurate prediction of the system solution.

\subsection{Important features of the parametrization}

\begin{figure}[ht!]
    \centering
    \includegraphics[width=1\columnwidth]{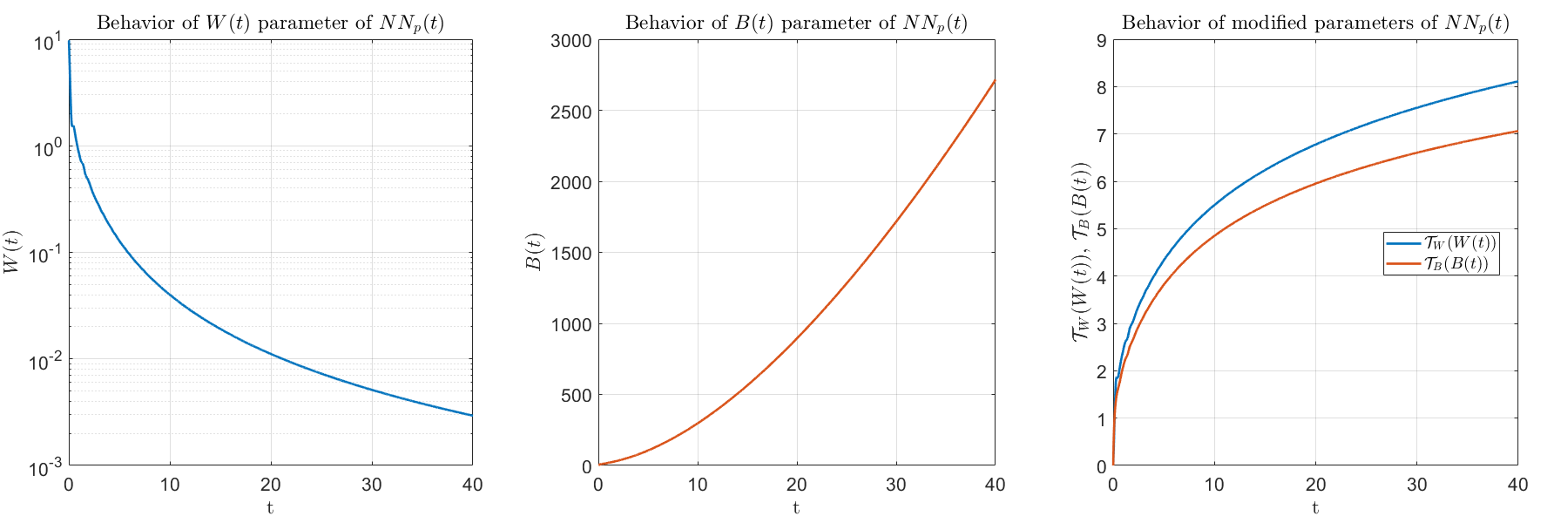}
    \caption{\label{fig:params_behavior} The behavior of $W(t)$ (left) and $B(t)$ (middle) parameters of ReLU-networks $NN_P(t)$, used to approximate the solution of Eqs. (\ref{problem:smol-bin}) with (\ref{monodisp}).  After applying the transformation \eqref{eq:params_transform}, these functions demonstrate more uniform and comparable behavior  (right).}
\end{figure}

Since $W(t)$ and $B(t)$ lack an explicit expression, the 
relevant information may be gained from the phenomenological insight. Taking a closer look at the behavior of $W(t)$ and $B(t)$ shown in Fig. \ref{fig:params_behavior} (left and middle panels) we notice that both functions are monotone and smooth. The curve shown in the left panel decays towards $0$ whereas the other curve appears to grow monotonically without reaching saturation. To make the behavior of these functions more 
uniform and comparable, we apply the following bijective transformations: 
\begin{align}
    \label{eq:params_transform}
    \mathcal{T}_W(W(t)) &= -(\log W(t) - \log W(0))\\
    \mathcal{T}_B(B(t)) &= \log B(t) - \log B(0) . \nonumber
\end{align}

As one can see from Fig. \ref{fig:params_behavior}, the functions, corresponding to the transformed variables $\mathcal{T}_W(W(t))$ and $\mathcal{T}_B(B(t))$ share an important property: they increase with time sub-linearly. A plausible conjecture, that they behave as $\sim t^{\gamma_{W/B}}$, with $\gamma_{W/B} <1$, implies alternating derivatives (a logarithmic dependence yields the same result). 
Namely, the first derivative of these functions is positive, the second -- negative, the third is again positive, and so on. 

\begin{align}
    \label{eq:params_modified_property}
    &\mathcal{T}_W(W(t))', \: \mathcal{T}_B(B(t))' > 0, \nonumber \\ 
    &\mathcal{T}_W(W(t))'', \: \mathcal{T}_B(B(t))'' < 0, \\
    &\mathcal{T}_W(W(t))^{(3)}, \: \mathcal{T}_B(B(t))^{(3)} > 0 \nonumber
\end{align}
This additional phenomenological knowledge, Eq. (\ref{eq:params_modified_property}), will be used to inform our predicting system about the behavior of the functions on the interval of interest $(\tau, T]$.

\subsection{Knowledge-informed extrapolation. Predictive neural network}

To construct the vector function ${\bf p}(t)= \left( W(t), B(t) \right)$, which could potentially extrapolate beyond the initial interval $[0,\tau]$ we need to find an appropriate class of systems (\ref{eq:parameters_system}). Ideally, such systems should allow efficient computation of the values of the vector ${\bf p}(t)$. At the same time, they should enable convenient and accurate fitting of ${\bf p}(t)$ to date over the interval $[0,\tau]$ for which the values of ${\bf p}(t)$ are known. 

In our current work, we utilize the fact that some classes of fee-dforward neural networks may be considered as closed-form solutions of (\ref{eq:parameters_system}), along with the corresponding ODE subsystem. Detailed justification of this fact is provided in the Appendix.

For the sake of simplicity, in what follows we will be focusing primarily on the reconstruction of the transformed parameters -- the functions $\mathcal{T}_W(W(t))$ and $\mathcal{T}_B(B(t))$ -- over the entire interval $[0,T]$. It turns out that the values of these parameters could be predicted by outputs of two independent/decoupled networks.  Given that the expected behavior of $\mathcal{T}_W(W(t))$ and $\mathcal{T}_B(B(t))$ is similar, the architectures and training protocols for these two networks were kept identical too. 
Below we provide a general overview of our training setup and results for the reconstruction of $\mathcal{T}_W(W(t))$, with all technical details provided in Appendix.

Let a neural network $NN_e$ with a scalar input $t$ and a scalar output $NN_e(t)$ be the network reconstructing the values of $\mathcal{T}_W(W(t))$. Network's outputs $NN_e(t)$ should satisfy constraints (\ref{eq:params_modified_property}) for $t \in [0, T]$ and must also reproduce the values of $\mathcal{T}_W(W(t))$ for  $t \in [0, \tau]$. To ensure that these requirements are met, we introduced two corresponding loss functions enforcing relevant sets of constraints. The first loss function, $L_d$, penalised the network for the deviation from the requirements listed in (\ref{eq:params_modified_property}) over the whole interval $t \in [0,T]$: 

\begin{equation}
    \label{eq:L_d_loss}
    L_d(NN_e) = \int_0^T \left(C_1  ReLU(-NN_e'(t)) + C_2  ReLU(NN_e''(t)) + C_3  ReLU(-NN_e^{(3)}(t)) \right)dt.
\end{equation}
In (\ref{eq:L_d_loss}),  $NN_e'(t)$, $NN_e''(t)$, and $NN_e^{(3)}(t)$ are the first, second and third derivatives of the output of the network, and $C_1,C_2,C_3\in\Real_{>0}$ are penalty coefficients.   

The second loss function, 
\begin{equation}
    \label{eq:L_u_loss}
    L_u(NN_e,\mathcal{T}_W,W) = \int_0^{\tau} |NN_e(t) - \mathcal{T}_W(W(t))|^2 d t, 
\end{equation}
penalised the deviation of the network's output from the values of $\mathcal{T}_W(W(t))$ on the interval $[0,\tau]$.

\begin{figure}[ht!]
    \centering
    \includegraphics[width=0.7\columnwidth]{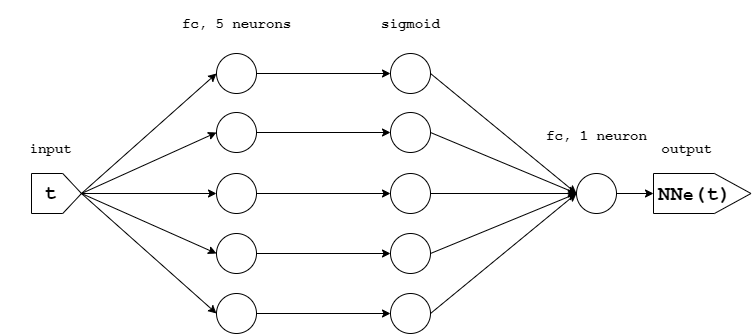}
    \caption{\label{fig:NN_e_structure} Structure of the predictive network $NN_e$ utilized for knowledge-informed extrapolation of the parameters.}
\end{figure}

In general, one can use any network (\ref{eq:NN_p_equation_gen}) to model $\mathcal{T}_W(W(t))$. Here we used a network $NN_e$ with two fully-connected layers with 5 neurons in the first layer, 1 neuron in the second, and sigmoid activation function between layers (see Fig. \ref{fig:NN_e_structure}). 

The choice of the sigmoid function is motivated by the fact that these functions can be viewed as solutions of an ordinary differential equation (see Appendix) and hence networks with such activation functions are capable of implementing systems of relevant equations -- the property which we exploit here. During network's training phase,  
we alternated the minimization of $L_u$ and $L_d$ at each epoch. As already noted, for the extrapolation of $\mathcal{T}_B(B(t))$, the same procedure has been exploited. 

\section{Knowledge-informed prediction of aggregation kinetics }\label{sec:prediction}

\subsection{Knowledge-informed neuro-solver for the kernel $K_{i,j}=1$}

To test the efficiency and the accuracy of the proposed approach numerically we considered system \eqref{problem:smol-bin}  with initial conditions \eqref{monodisp} and reaction kernel \eqref{Kij1}. Results of this numerical exploration are provided below, and the code implementing the experiments can be accessed at \cite{know-info-code}.

To construct extrapolations of $\mathcal{T}_W(W(t))$ and $\mathcal{T}_B(B(t))$ we used the values of solutions of \eqref{problem:smol-bin} computed over the interval $t \in [10, 19]$. The values from the interval $t \in [0, 10]$ were omitted due to observed non-monotonicity of $\mathcal{T}_W(W(t))$ and $\mathcal{T}_B(B(t))$ on this initial ``transien'' interval. The values from the interval $t \in [19, 20]$ were used for the validation of the quality of prediction during training. 

\begin{figure}[ht!]
    \centering
    \includegraphics[width=0.7\columnwidth]{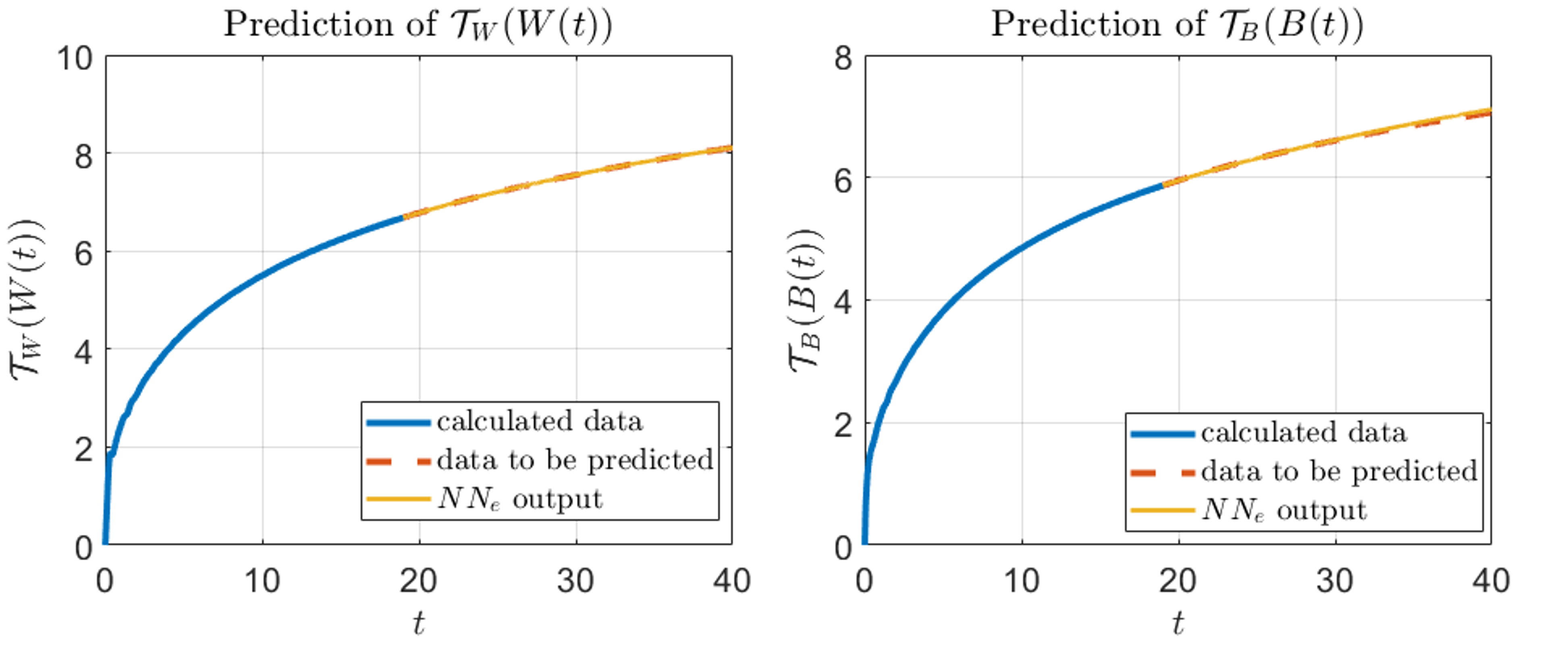}
    \caption{\label{fig:params_predicted} Extrapolation of modified parameters $\mathcal{T}_W(W(t))$ (left) and $\mathcal{T}_B(B(t))$ (right) made with the network $NN_e$.}
\end{figure}

After completing the training, we applied the predictive network to extrapolate the values of $\mathcal{T}_W(W(t))$ and $\mathcal{T}_B(B(t))$ over $t\in [20,40]$. The results of extrapolation are depicted in Fig. \ref{fig:params_predicted}. Notice that the size of the target interval, $[10,40]$, is three times larger than the size of the interval on which the network was trained.  

The values of $W(t)$ and $B(t)$ were restored from the values of $\mathcal{T}_W(t)$, $\mathcal{T}_B(t)$ in accordance with the following inverse transformation:
\begin{equation}
    W(t) = W(0) \cdot e^{-\mathcal{T}_W(t)}, \qquad \qquad B(t) = B(0) \cdot e^{\mathcal{T}_B(t)}.
    \end{equation}

\begin{figure}[ht!]
    \centering
    \includegraphics[width=1\columnwidth]{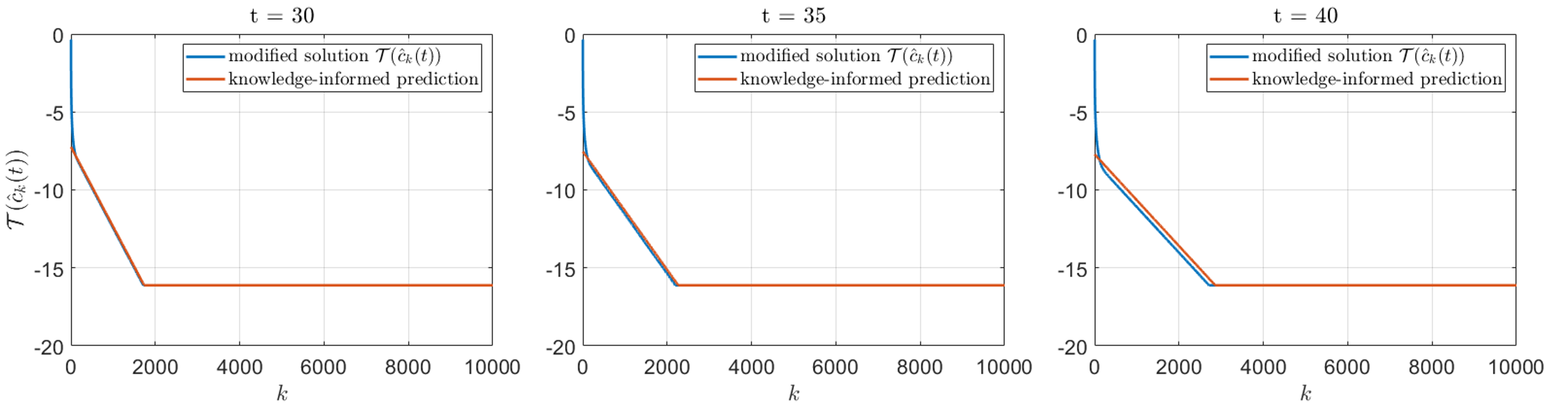}
    \caption{\label{fig:NN_p_predict} Knowledge-informed prediction of the modified solution $\mathcal{T}(\hat{c}_k(t))$ made with the network $NN_p$ with predicted free parameters $W(t)$ and $B(t)$. The training time interval is $[10,20]$.}
\end{figure}

This, together with the network $NN_p$, results in the predicted solution of the aggregation equations (\ref{problem:smol-bin}) over the target interval $[10,40]$. The values of $\mathcal{T}(\hat{c}
_k(t))$ and the predicted solution are shown in Fig. \ref{fig:NN_p_predict}. The figure indicates a fairly accurate prediction of the cluster densities over time, demonstrating the efficiency of the method.

The above analysis has been performed for the time interval, associated with a moderate computational cost. To study the predictive capacity of knowledge-informed neuro-solver for longer times, we expand the prediction window for the same kernel, $K_{i,j}=1$,  from $t \in [20, 40]$ to $t \in [20, 80]$.  This is on the boundary of computational power and memory available in ordinary computers. We used the same setup as above: the training set consisted of the values of parameters $W(t)$ and $B(t)$ of the parametrizing network $NN_p$ for $t \in [10, 20]$ and  the target interval was $t \in [10, 80]$. That is, the target interval was 7 times larger than the training one. 

\begin{figure}[ht!]
    \centering
    \includegraphics[width=0.7\columnwidth]{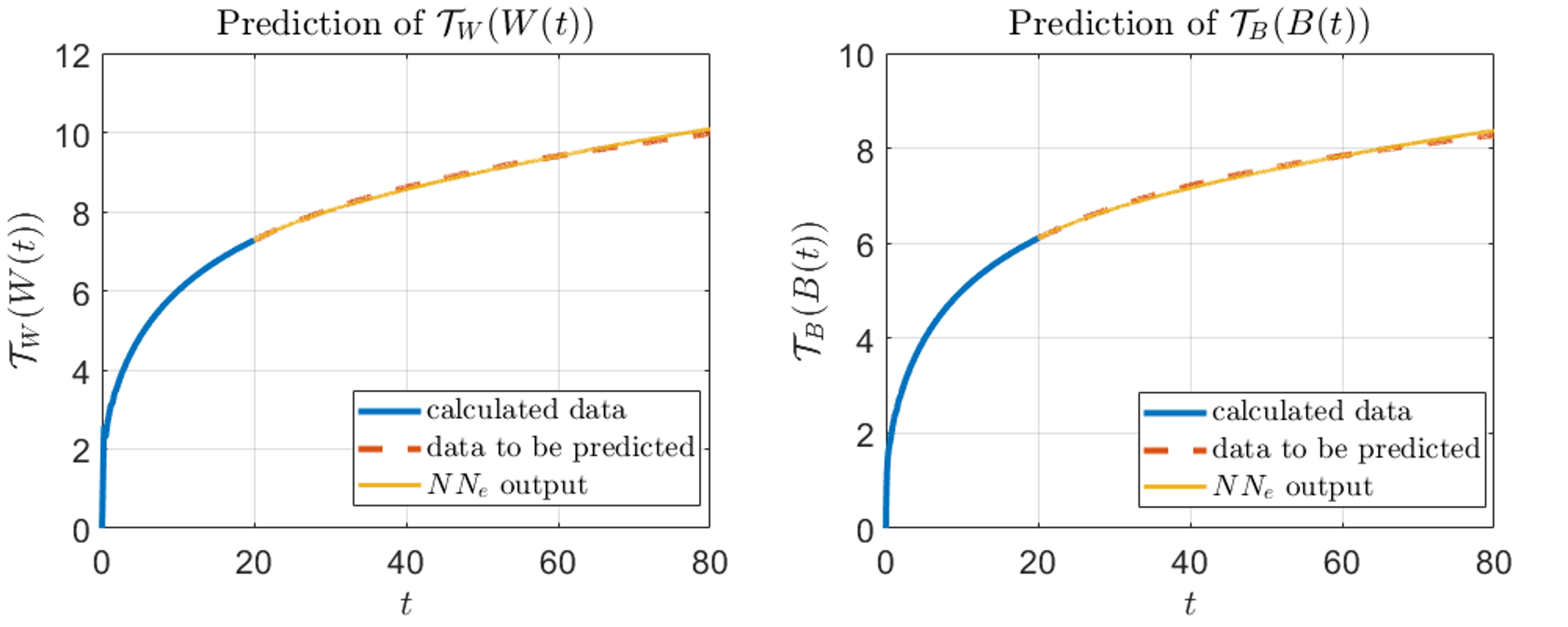}
    \caption{\label{fig:params_predicted_cont}  Extrapolation of modified parameters $\mathcal{T}_W(W(t))$ (left) and $\mathcal{T}_B(B(t))$ (right) for aggregation process with $K_{i,j} = 1$ on expanded time interval with $t \in [0, 80]$.}
\end{figure}

\begin{figure}[ht!]
    \centering
    \includegraphics[width=1\columnwidth]{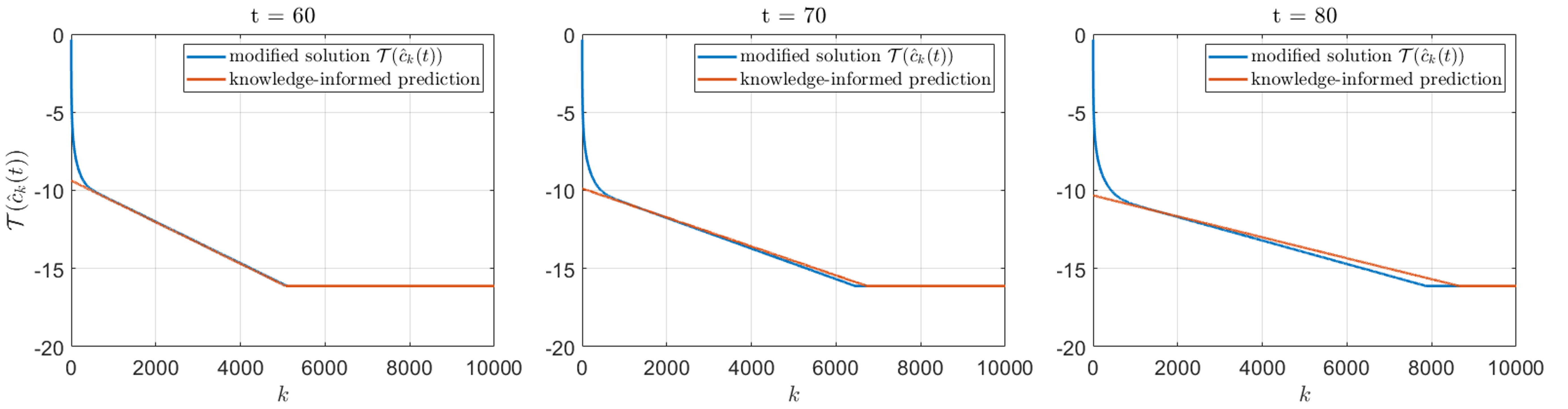}
    \caption{\label{fig:NN_p_predict_cont} Knowledge-informed prediction of the modified solution $\mathcal{T}(\hat{c}_k(t))$ for the aggregation process with unit kernel $K_{i,j} = 1$. The training time interval is $[10,20]$.}
\end{figure}

Fig. \ref{fig:params_predicted_cont} contains comparison of actual values of parameters $W(t)$ and $B(t)$ of parametrizing network and values extrapolated from the interval of $t \in [10, 20]$ to the interval of $t \in [10, 80]$. The results of extrapolation are relatively accurate, which leads to the accurate prediction of the distribution of particles $\mathcal{T}(\hat{c}_k(t))$ reconstructed from the parameters. Comparison of predicted solution and the actual solution is depicted in Fig. \ref{fig:NN_p_predict_cont}.

Note that the exploration of the system behavior, for such a large maximum $t$, is associated with a dramatic increase of the computational cost. It increases much faster than linear with increasing time, as one needs to account for larger and larger emerging particles, due to increasing truncation error. In the case of maximal time equal $80$ we have to consider particles of the size up to $40\, 000$. Hence, we need to solve a system of more than  $40\, 000$ equations, which is computationally challenging, both in terms of the computation time and requested memory.

\subsection{Knowledge-informed neuro-solver for other kernels}

\begin{figure}[ht!]
    \centering
    \includegraphics[width=0.7\columnwidth]{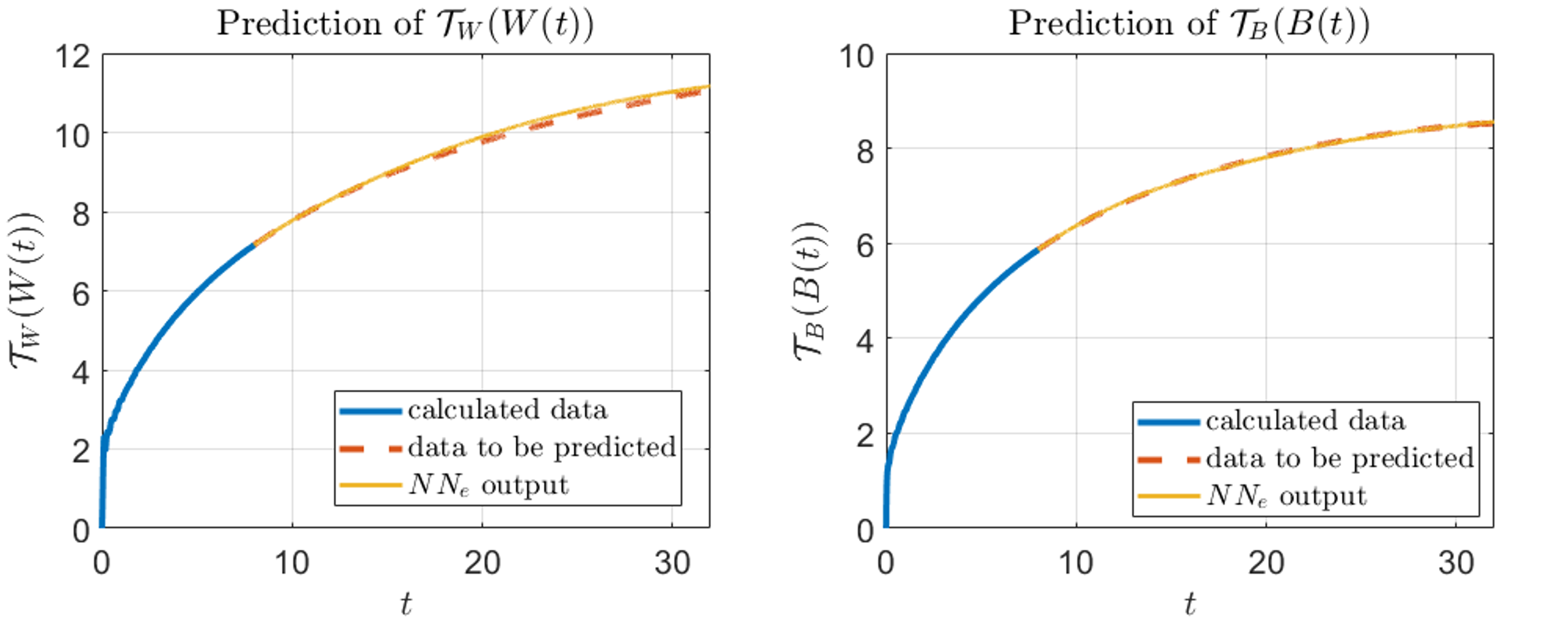}
    \caption{\label{fig:prod_ker_params_predicted} Extrapolation of modified parameters $\mathcal{T}_W(W(t))$ (left) and $\mathcal{T}_B(B(t))$ (right) for the aggregation process with the product aggregation kernel $K_{i,j} = (ij)^{0.2}$.}
\end{figure}

\begin{figure}[ht!]
    \centering
    \includegraphics[width=1\columnwidth]{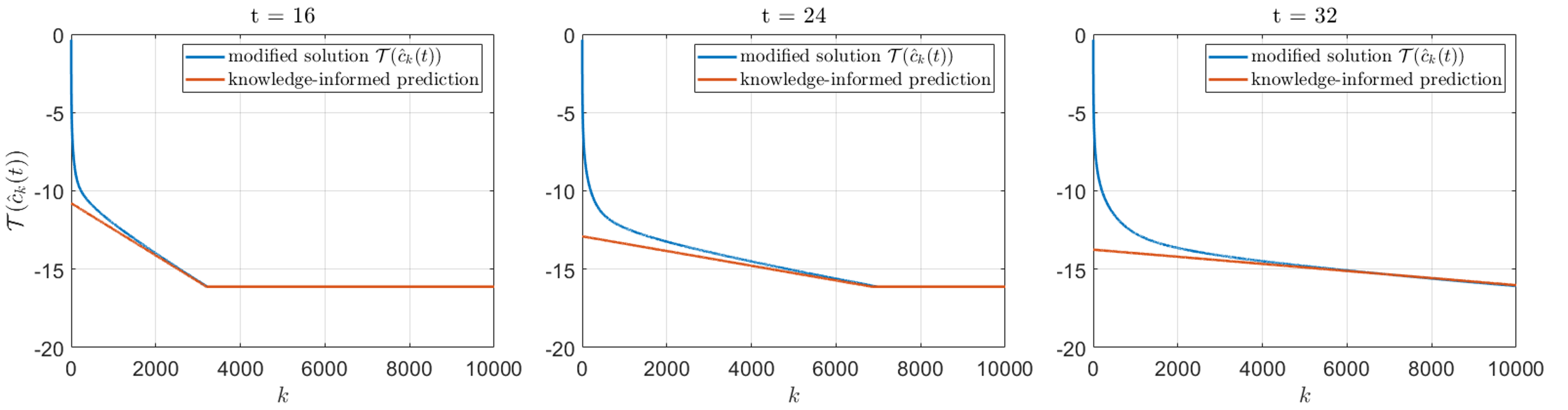}
    \caption{\label{fig:NN_p_prod_ker_predict} Knowledge-informed prediction of the modified solution $\mathcal{T}(\hat{c}_k(t))$ for the aggregation process with the product aggregation kernel $K_{i,j} = (ij)^{0.2}$. The training time interval is $[4,8]$.}
\end{figure}

To explore the applicability of the proposed method to other kernels we considered the solution of the aggregation equations \eqref{problem:smol-bin} with the same initial conditions \eqref{monodisp}, but with the reaction kernels from the class \eqref{eq:Cijgen} with $\nu \neq 0$ and $\mu \neq 0$. Namely, we use the product kernel, with $\nu = \mu = 0.2$:

\[
K_{i,j} = (ij)^{0.2}.
\]
As mentioned in Section \ref{sec:formal}, kernels of this class are widely used in studies of aggregation processes and phenomena. The results of the application of the method to equations with the product kernel are shown in Fig. \ref{fig:prod_ker_params_predicted} - \ref{fig:NN_p_prod_ker_predict}. 

Note that, since the aggregation occurs significantly faster in systems with this kernel, the time scale is different. To further explore computational efficiency, similar to the previous case, the target interval of the prediction $[4, 32]$ is seven times larger than the size of the training interval $[4,8]$. Note that in this case, computational cost increases even faster, as kernels from the considered family have more complex low-rank approximations. Again we see a good prediction of parameters, albeit the prediction of solutions is limited by the simplicity of the parameterizing network.

\begin{figure}[ht!]
    \centering
    \includegraphics[width=0.7\columnwidth]{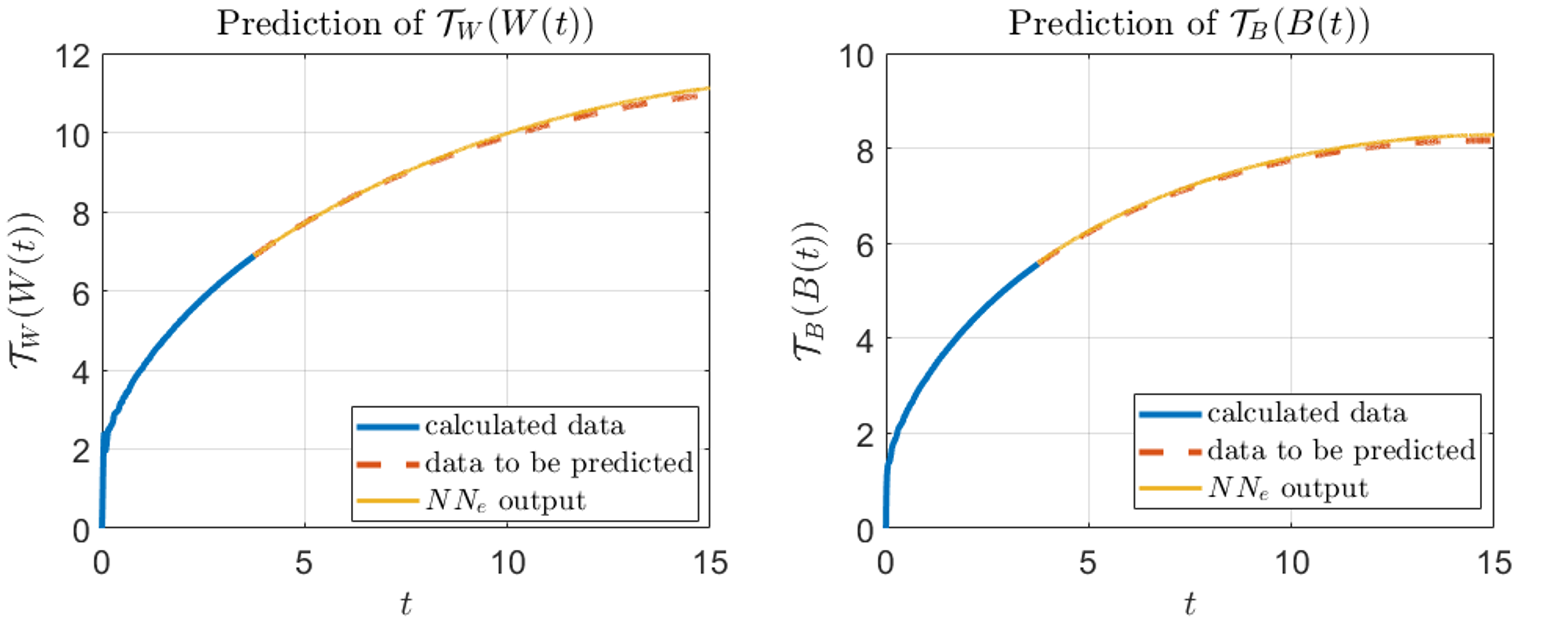}
    \caption{\label{fig:sum_ker_params_predicted} Extrapolation of modified parameters $\mathcal{T}_W(W(t))$ (left) and $\mathcal{T}_B(B(t))$ (right) for the aggregation process with the sum-like aggregation kernel $K_{i,j} = i^{0.5} + j^{0.5}$.}
\end{figure}

\begin{figure}[ht!]
    \centering
    \includegraphics[width=1\columnwidth]{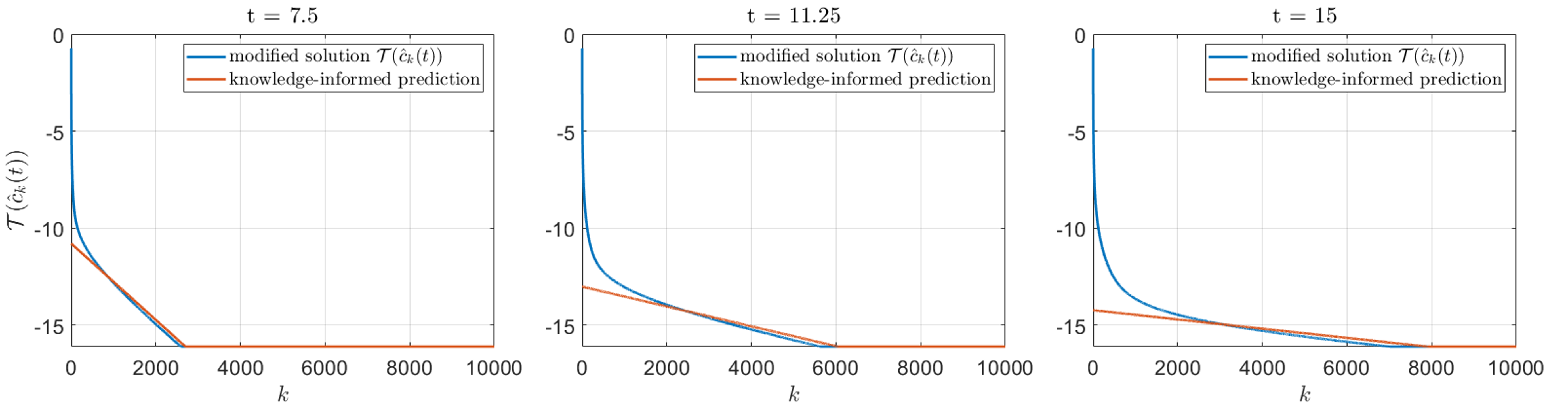}
    \caption{\label{fig:NN_p_sum_ker_predict} Knowledge-informed prediction of the modified solution $\mathcal{T}(\hat{c}_k(t))$ for the aggregation process with the sum-like aggregation kernel $K_{i,j} = i^{0.5} + j^{0.5}$. The training time interval is $[1.875,3.75]$.}
\end{figure}

Next we applied our approach to the summation kernel, with $\nu = 0.5$ and $ \mu = 0$:
\[
K_{i,j} = i^{0.5} + j^{0.5}.
\]
In this case, the training interval was $[1.875, 3.75]$, and the target interval was $[1.875, 15]$. The results of the knowledge-informed prediction are illustrated in Fig. \ref{fig:sum_ker_params_predicted}, \ref{fig:NN_p_sum_ker_predict}. The predictions again are initially accurate with slight divergence from the true solutions towards the end of the target interval. 

\subsection{Knowledge-informed neuro-solver of higher order}

\begin{figure}[ht!]
    \centering
    \includegraphics[width=0.5\columnwidth]{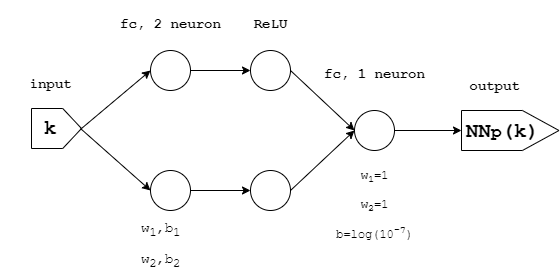}
    \caption{\label{fig:big_NN_p_structure} Structure of parametrizing network $NN_p$ utilized for more accurate prediction of both distributions of particles with lower and higher masses.}
\end{figure}

To further explore the potential of the proposed technique we built a knowledge-informed solver utilizing different parametric family $F(k, \mathbf{p}(t))$, i.e. parametrizing network $NN_p$. As one might see from Fig. \ref{fig:NN_p_predict}, simple one-neuron parametrizing network fails to accurately predict behavior of particles with small masses. To improve the predictive ability, we will use two-neuron neural network with the structure depicted in Fig. \ref{fig:big_NN_p_structure}.

Essentially, a new parametrizing network is a parametric family of partially-linear functions with three linear segments. As the new family includes previously utilized one, with only two segments, we expect increase of predictive capability.

Since a noticeable difference, between the distribution of particles of small and large mass started to appear quite late, at large times, we used a different dataset for the new parametrizing network. For the model with kernel $K_i,j = 1$ we utilized, the size distributions $\hat{c}_k(t)$ for $t$ from 50 to 80. The training set used the data for $t$ from 50 to 60, while the knowledge-informed solver made the prediction for $t$ from 60 to 80.

\begin{figure}[ht!]
    \centering
    \includegraphics[width=1\columnwidth]{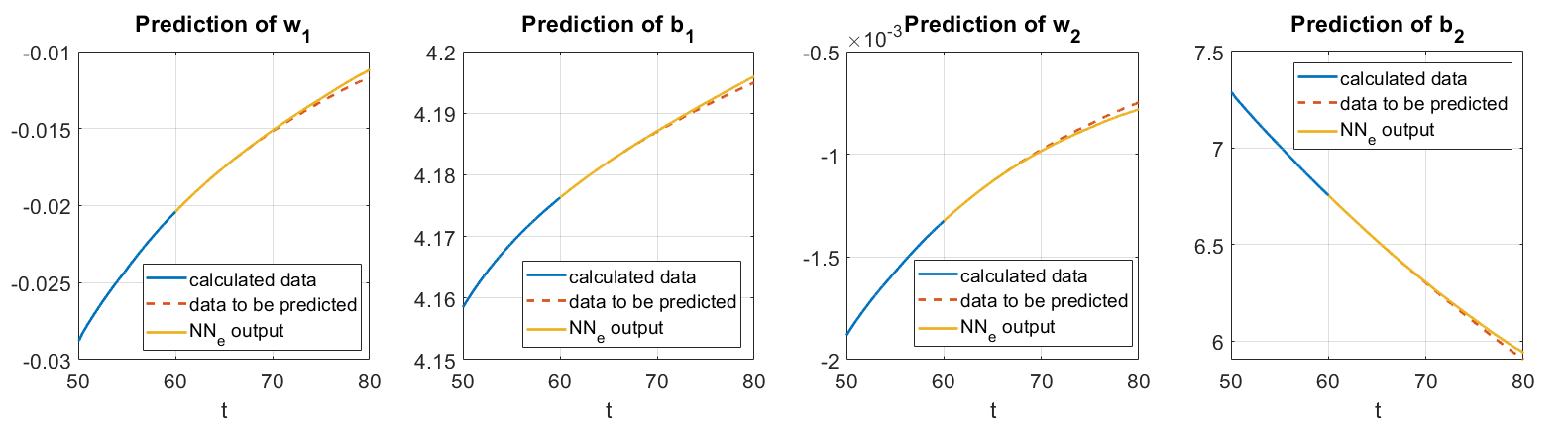}
    \caption{\label{fig:big_net_params_predicted} Extrapolation of parameters of more complicated parametrizing network $NN_p$.} Note, that in contrast to other plots of similar quantities, we do not use here the logarithmic scale for $w_{1/2}$ and $b_{1/2}$.
\end{figure}

\begin{figure}[ht!]
    \centering
    \includegraphics[width=1\columnwidth]{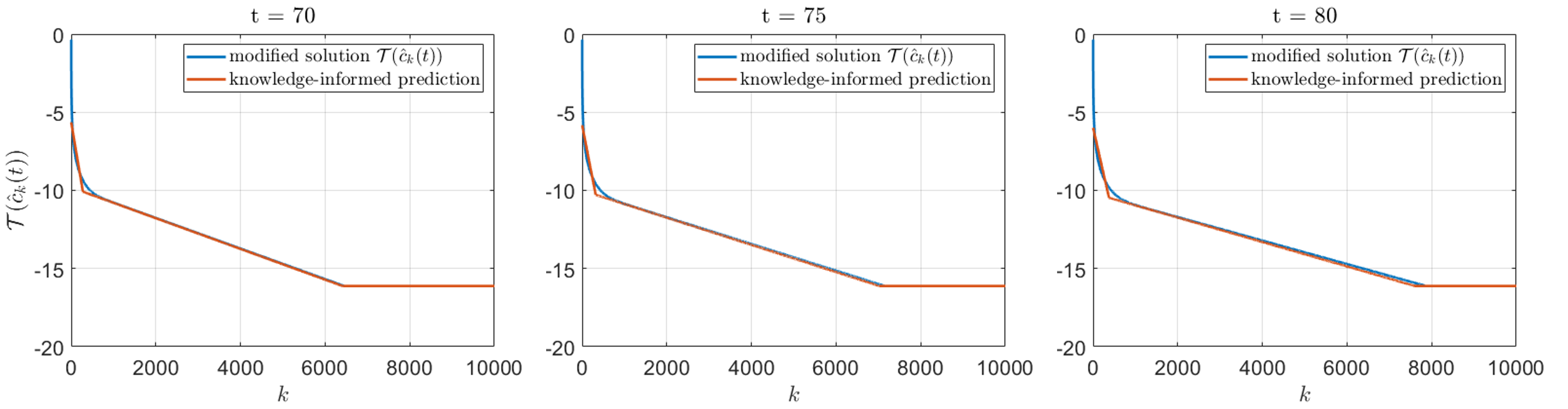}
    \caption{\label{fig:big_net_predicted} Knowledge-informed prediction of the modified solution $\mathcal{T}(\hat{c}_k(t))$}
\end{figure}

The prediction of four learnable parameters is depicted in Fig. \ref{fig:big_net_params_predicted}.  As one can  see from Fig. \ref{fig:big_net_predicted}, the proposed more complex parametrizing network, provides a substantial improvement for  the prediction of the small-mass particle behavior. 

\subsection{Computational efficiency and error estimates of knowledge-informed prediction}

A theoretical bound on the computational costs of the proposed method can be constructed from the following three additive components. The first component is the cost of computing solutions over the initial interval $[0,\tau]$ using state-of-the-art conventional solvers such as \cite{matveev2015}. These costs scale as $O((N\log N)N_\tau)$ \cite{matveev2015}, where $N$ is the system's size and $N_\tau$ is the number of temporal steps in the interval $[0,\tau]$. The second component is the overall computational cost of building parameterizing and predicting networks. Since we used a standard gradient algorithm to determine the weights and biases of the networks, these costs are proportional to the number of operations needed to compute gradients of loss functions (\ref{eq:L_d_loss}), (\ref{eq:L_u_loss}). These costs increase no faster than $O(m N_T)$, where $m$ is the sum of the number of parameters in both networks, and $N_T$ is the number of temporal steps in the interval $[0,T]$. The third component is the cost of computing an approximation of the solution over $[0,T]$. Since the approximation is computed by the above neural networks, this cost grows no faster than $O(m N_T)$. Hence the total costs are at most
\[
\underbrace{O((N\log N)N_\tau)}_{\mathrm{Precalculation}} + \underbrace{O(m N_T)}_{\mathrm{Retrieving}+\mathrm{Prediction}}.
\]

To evaluate the actual computational efficiency of our approach, we computed and compared its CPU time with that of the state-of-the-art methods. When running the experiments, we did not involve GPUs or used dedicated multi-thread optimisation for parameter retrieving and solution prediction in the proposed knowledge-informed solver. We have explored three different systems with kernels $K_{i,j} = 1$, $K_{i,j} = (i \cdot j)^{0.2}$,  and  $K_{i,j} = i^{0.5} + j^{0.5}$. These systems modeled aggregation of particles with sizes up to $40 \,000$, using $40 \,000$ equations. The time horizon of modeling was limited to $8000$ time steps.

For the computation of "true" solution for these three cases we used fast-solvers, based on low-rank approximation of the reaction rate kernels, described in \cite{matveev2014fasttransliteration, matveev2015fast}. Computation time for the first, second and third case were respectively,  $480.20s$, $1981.40s$, and $2054.20s$ (see also Table \ref{tab:computation-times}). Note,  that the implementation of the fast-solver is parallelized to processor cores, while the implementation of our approach is fully sequential. This implies the potential for further acceleration of our method.

The application of the knowledge-informed approach requires preliminary numerical integration. In all three cases, we used a state-of-the-art fast solver to integrate the model over the first $2000$ time steps. Usage of only $2000$ time steps instead of $8000$, the modelling horizon in our examples, leads to dramatic reduction of computational costs. The error related to truncation \eqref{eq:Smol-trunc} of the initially infinite system \eqref{eq:Smol1} increases with time, hence we can utilize much smaller system to compute the "true" solution. The system used to compute the first $2000$ time steps contained $10 \,000$ equations, four times smaller than the one needed to compute $8000$ steps.
For the three cases, the corresponding computation times for these pre-calculations were, respectively, $22.40s$, $94.30s$ and $92.70$s. 

Knowledge-informed prediction consists of two main parts: retrieving of parametrizing network parameters, and prediction of these parameters for the target interval. For all three cases, the corresponding computation time was essentially the same. Parameter retrieving was done on the same number of 1000 time steps and took $33s$, with the variation of less than $2s$. Parameter prediction involving two independent runs of training for two parameters took, on average, $252s$, with the variation of less than $33s$. 
The information about the CPU time used for different computations above is presented in Table \ref{tab:computation-times}.
\begin{table}[]
    \centering
    \begin{tabular}{| c | c | c | c | c | c |}
    \hline
    & Fast-solver & Knowledge-Informed & Precalculation & Retrieving & Prediction \\
    \hline
    $K_{i,j}=1$, $N=40000$, $T_{st} = 8000$ & $480.20s$ & $276.10s$ & $22.40s$ & $32.90s$ & $220.80s$ \\
    \hline
    $K_{i,j}=(i \cdot j)^{0.2}$, $N=40000$, $T_{st} = 8000$ & $1981.40s$ & $377.60s$ & $94.30s$ & $31.70s$ & $251.60s$ \\ 
    \hline
    $K_{i,j}=i^{0.5} + j^{0.5}$, $N=40000$, $T_{st} = 8000$ & $2054.20s$ & $411.10s$ & $92.70s$ & $33.40s$ & $285.00s$ \\ 
    \hline
    \end{tabular}
    \caption{\label{tab:computation-times} Comparison of CPU time spent on calculation of aggregation processes with fast-solver and knowledge-informed neuro-integrator. Computation time of knowledge-informed approach is divided into three parts: time spent for pre-calculation of the training data, time spent on retrieving parameters of the parametrizing network, and time spent for the prediction of these parameters.}
\end{table}

As one can see from the Table, for the simple kernel $K_{i,j} = 1$ the knowledge-informed approach outperforms the conventional solver and reduces, almost twice, the computational time. For systems with more complex kernels, the knowledge-informed approach provides the solution almost five times faster than the conventional fast-solver.

The advantage of the knowledge-informed approach over the traditional fast solvers becomes even more pronounced when the modelling horizon $T$ grows large. This is because large modeling horizons require larger sizes of the truncated systems to retain desired levels of modelling accuracy. This, in turn, implies superlinear growth of the computational costs of the fast solver.  On the contrary, the knowledge-informed approach does not require simulations over large intervals of time. It can, in principle, utilize information from solutions of the smaller-sized truncated systems.  The only part in the proposed approach that depends on the size of the target interval (modelling horizon $T$) is the training of the extrapolating network. This dependence is through the computation of gradients of the loss function (\ref{eq:L_d_loss}) over the entire interval $[0,T]$. The computational costs of the latter (very cheap operations) grow linearly with $T$. 

The knowledge-informed prediction, utilizing more complex network has been included to only illustrate the applicability of the proposed method to other parametric families. The analysis of its efficiency in beyond the main scope of the article. As more complex behavior of the solution requiring more complex parametrizing network, the preliminary calculations take more computation time.  For the considered case with the kernel $K_{i,j} = 1$,  the pre-calculations made with $30 \,000$ equations, take 240.40s, the parameter retrieving takes 31.60s, and parameter prediction takes  410.20s. This, in total, exceeds the CPU times of 480.20s of the conventional fast-solver. This time, however, corresponds to fully sequential computing, while the fast solver computations were parallelized. The parallelization of the parameter prediction stage will dramatically decrease the computation time, leading to the superiority of the new approach also for high-order neuro-solver. 

Theoretical background supporting the knowledge-informed approach, including a-priory error estimates, is provided in Appendix, Sections A and B. In Section A  we recall that neural networks with sigmoid functions are universal approximators and also can be modelled by equivalent systems of ordinary differential equations (Theorem \ref{thm:neural_network:dynamic}). In Section B we present a general extrapolation result capturing the requirements and assumptions needed to achieve successful extrapolation (Theorem \ref{thm:exptrapolation:general}). This result, together with Theorem \ref{thm:neural_network:dynamic}, is used to produce extrapolation bounds for neural networks (Corollary \ref{cor:neural_net_error}).
Although these estimates are largely qualitative, the general expressions allow us to determine the main factors contributing to controlling the extrapolation error. 

\section{Discussion}\label{sec:discussion}

A novel approach for modelling aggregation kinetics, described in this work, showed promising results in a view of the balance between the achieved accuracy of prediction and the overall structural simplicity of the reduced model. We note that the  application of the method in its simplest form is not very accurate for particles of relatively small size, see see Fig. \ref{fig:NN_p_predict}, \ref{fig:NN_p_prod_ker_predict} and  \ref{fig:NN_p_sum_ker_predict} for  small $k$. This however, is not so important. Indeed, our main goal is predicting densities of  large-size particles. Densities of small particles are easily accessible for short times, by standard techniques, applied for a small set of ODE. For sufficiently large times, only a negligible fraction of total mass is associated with small particles. Furthermore, one can consider an extended class of functions $\mathcal{F}$, e.g., by including several ReLU nonlinearities. Using more complex parameterizing networks, for instance, networks generating a family of piece-wise linear or even polynomial functions, will produce a better result,  albeit at higher computational costs. We have demonstrate this by utilizing the knowledge-informed solver of high order. Applying the two-neuron, instead of one-neuron network,  for the parametrizing ANN,  resulted in a good  accuracy for the particles densities of all sizes. 

In the current study we illustrated the application of the method with examples when the modelling horizon $T$ was only several times larger than the size of the ``training'' time interval $[0,\tau]$. This was due to the computational constraints of the comparison protocol: to evaluate the difference between solutions produced by our method and the standard solver one needs to apply the standard solver to the entire interval $[0,T]$ first. This is a major computational task for a standalone machine. The new method does not have these limitations. Nevertheless, since we could not evaluate the difference between the standard and the new method, we did not show predictions for $T \gg \tau$.

The application of the general framework of the method described in Section \ref{sec:formal} is not constrained by the type of the aggregation kernel or the form of the aggregation equations. This is in sharp contrast with the fast-solver, applicable only to kernels reducible to a sum of few low-rank components. Other classes of kinetic equations can be also explored. Moreover, one can use other underlying dynamical systems, in addition, or instead of the logistic equations, utilized here. This may enable an efficient application of knowledge-informed prediction, which opens a potentially vast field of fast neuro-integrators, capable of modelling other  important for application processes.

We would also like to mention that the proposed neuro-integrator method bears a degree of similarity with a plethora of conventional numerical methods, such as the finite elements or the method of collocations in that it uses a piece-wise approximation of the original solutions. The difference, however, is that the traditional methods exploit a prior knowledge about differential equations to generate numerical solutions. Our method does not require such information. It exploits the information from the numerical values of the solution over short intervals of time, regardless of how these solutions were obtained. Other related methods are integration methods with adaptive moving meshes (see \cite{budd2009adaptivity} for a comprehensive review). In these methods, the temporal evolution of the mesh -- an analogue of parameters of the scaling functions in our case, is guided by the requirements to predict the location of ``interesting'' behaviors over several steps ahead. As the result, the mesh concentrates around the most relevant points in space and time,  ensuring minimal errors. Here, we focus on the possibility of predicting solutions over long intervals of time and without the need to know the equations themselves. If the equations are known, then the proposed method could be used as a part of the schemes of ``large steps'', see chapter 9.1 of  \cite{gorban2005invariant} and 
Refs. \cite{gorban2004uniqueness,gorban1996relaxational}. An approximate solution at a large-size time step is produced there, followed by a projection of the latter onto a relevant subspace of true solutions. Exploring this direction, however, is beyond the scope of the current work. 

\section{Conclusion}\label{sec:conclusion}

In this work, we developed a novel approach for modelling aggregation kinetics, for time intervals inaccessible for the existing  computational methods. The main idea of the method is based on a specific ("smart") application of artificial neural networks (ANN). These extrapolate the numerical solution of a relatively small ODE system, describing agglomeration for a limited time interval, to much larger systems and significantly longer time.  Specifically, in our study, we demonstrated the feasibility of the method when the modelling horizon was four times larger than the initial, training interval. We did not evaluate the method over longer modelling horizons due to the computational constraints stemming from the application of the standard solver to larger intervals.

To enable extrapolation beyond the initial interval, we supply the standard ANN training with additional information about expected qualitative behavior of the solution outside the training interval. This motivates the name, "knowledge-informed", of our approach. In practice, we deploy two entangled ANNs -- the parametrizing and predictive network. The former network is used for an appropriate parametrization of the numerical solution of a relatively small system  of $N$ ODEs, for a short time interval, within  a chosen class of functions;  here --  the feedforward neural networks. 
At this step we, essentially,  make a dimensional reduction of the problem. Indeed, the solution vector ${\bf C}(t)$ of the ODE, with $N$ components, associated with the densities of aggregates of size $1, \, 2, \ldots N$, is mapped on the vector of parameters, ${\bf p}(t)$, of much smaller dimension. Knowing ${\bf p}(t)$ and the class of parametrizing functions, one reconstructs ${\bf C}(t)$. Then, the predictive  network performs the extrapolation of the vector ${\bf p}(t)$ for much longer (target) time interval, utilizing the knowledge about the time dependence of  ${\bf p}(t)$, extracted from its behavior on the short time interval. Here we can use very small networks -- only two neurons for the  parametrizing ANN and six for the predictive ANN. We apply our new method to model aggregation kinetics for a couple of reactive kernels, including the Smoluchowski kernel, $K_{i}={\rm const.}$, and other kernels, important for applications -- the product kernel, $K_{ij}=(ij)^{\mu}$, and the summation one, $K_{ij}=i^{\nu}+j^{\nu}$ (we use $\mu=0.2$ and $\nu=0.5$). For all these kernels we obtain, with the new neuro-solver, a fairly accurate solution of the aggregation equations.

Surprisingly, very small, entangled ANN provide rather accurate prediction of the aggregation kinetics for relatively long time. With the use of the new neuro-solver, we obtain a solution of a large system of aggregation ODE for a large period of time, based on the solution of much smaller ODE system for shorter time. The utilized dimensional reduction allows application of ANNs, with a small number of neurons, and, eventually, fast computations associated with the neural networks. The reported method significantly outperforms the most efficient conventional state-of-the-art solver based on low-rank decomposition. The computational efficiency of the new method is especially important for solving aggregation equations for which the low-rank decomposition of the reaction kernels $K_{ij}$  fails (the ballistic kernel is a prominent example, see, e.g. \cite{Krapivsky,Leyvraz}). In this case one encounters the problem of $N^2$-complexity and the lack of available methods to overcome it; we expect that the new method will enable to address this problem.

\begin{acknowledgments}

The authors acknowledge the support by the Russian Science Foundation grant No.~21-11-00363, https://rscf.ru/project/21-11-00363/. Ivan Tyukin was supported by the UKRI Turing AI Fellowship EP/V025295/2.

\end{acknowledgments}

\bibliography{KnowInfo_Smol}
\bibliographystyle{ieeetr}

\section*{Appendix}

In this section, we provide relevant technical details, assumptions, and a-priory extrapolation error bounds for the analysis of the asymptotic properties of the proposed method. We also provide a detailed description of the training protocol and architectures of neural networks used in experiments.

Throughout the section, symbol $\Real$ denotes the field of real numbers, $\Real_{\geq 0}=\{x\in\Real| \ x\geq 0\}$, $\Real_{> 0}=\{x\in\Real| \ x> 0\}$, and $\Real^n$ denotes the $n$-dimensional real vector space, $\Natural$ denotes the set of natural numbers; $({x},{y})=\sum_{k=1}^n x_{k} y_{k}$ is the inner product of ${x}\in\Real^n$ and ${y}\in\Real^n$, and $\|{x}\|=\sqrt{({x},{x})}$ is the standard Euclidean norm  in $\Real^n$. By $\mathcal{C}^0([a,b],\Real^n)$,  $a,b\in\Real$, $a\leq b$, we will denote the space of continuous functions defined on the interval $[a,b]$ and taking values in $\Real^n$, $n\in\Natural$. $\mathcal{C}^p([a,b],\Real^n), \ p\in\Natural$, is the space of functions from $[a,b]$ to $\Real^n$ which have at least $p$ continuous derivatives. Let $\psi\in\mathcal{C}^0([a, b], \Real^n)$, symbol $\|\psi\|_{\infty,[a,b]}$ stands for $\max_{t\in[a,b]}\|\psi(t)\|$.

\subsection{Neural networks as systems of ordinary differential equations}\label{sec:nerual_networks_as_dynamical_systems}

Consider the following initial value problem
\[
\frac{dx}{dt}=w x (1-x), \ x(0)= x_0, \ w\in\Real, \ x_0\in\Real, \ x_0\in(0,1).
\]
The problem has a unique solution $x(t)=(1+\exp(-(w t + b)))^{-1}$ with
\[
b=- \ln \left(\frac{1}{x_0} - 1\right).
\]
It is hence clear that (cf. \cite{tyukin2003parameter}) that the output of the 2-layer network 
\begin{equation}\label{eq:reference_net}
NN_p=L^2 \sigma L^1:  \ NN_p(t)=\sum_{j=1}^{N_1} W^2_{j} \sigma(W^1_j t + b^1_j) + b^2,
\end{equation}
where $\sigma(s)=(1+\exp(-s))^{-1}$, coincides with 
\begin{equation}\label{eq:equivalent_odes_1}
y(t, W^2, W^1, b^2, b^1)=\sum_{j=1}^{N_1} W^2_j x_{j}(t; W^1,b^1) + b^2
\end{equation}
in which the functions of $x_j(\cdot;W^1,b^1)$, $j=1,\dots,N_1$ are unique solutions of
\begin{equation}\label{eq:equivalent_odes_2}
\begin{split}
&\left\{\begin{array}{ll}
\frac{dx_1}{dt}&= W^1_{1} x_1 (1-x_1)\\
&\cdots\\
\frac{dx_{N_1}}{dt}&= W^1_{N_1} x_{N_1} (1-x_{N_1})
\end{array} \right.\\
& \\
& x_j(0)\in (0,1) \ \mbox{for all} \ j=1,\dots,N_1
\end{split}
\end{equation}
with
\begin{equation}\label{eq:equivalent_odes_3}
b^1_j=-\ln \left(\frac{1}{x_j(0)}-1\right), \  j=1,\dots,N_1.
\end{equation}

In view of (\ref{eq:equivalent_odes_1})--(\ref{eq:equivalent_odes_3}) and \cite{cybenko1989approximation}, one can immediately state the following universal approximation theorem for the class of functions described by (\ref{eq:equivalent_odes_1}).

\begin{thm}{(See \cite{tyukin2003parameter})}\label{thm:neural_network:dynamic} Let $T>0$ be any positive real number. Functions (\ref{eq:equivalent_odes_1}) are dense in the space $\mathcal{C}^0([0,T],\Real)$ of real-valued continuous functions defined on $[0,T]$. That is for any $f:[0,T]\rightarrow \Real$, $f\in \mathcal{C}([0,T],\Real)$ and any $\varepsilon>0$ there exist $N_1\in\Natural, W^1\in\Real^{N_1},W^2\in\Real^{N_1}, b^2\in\Real, b^1\in\Real^{N_1}$ such that
\[
|f(t)-y(t, W^2, W^1, b^2, b^1)|<\varepsilon \ \mbox{for all} \ t\in[0,T].
\] 
\end{thm}

Equation (\ref{eq:reference_net}) is a closed-form expression of (\ref{eq:equivalent_odes_1}) -- (\ref{eq:equivalent_odes_3}) and hence shares all properties of the latter that are relevant for the task of extrapolation. These include continuous dependence of solutions on parameters and initial conditions and the unique dependence of the values of $x_j(t; W^1, b^1)$, $t\in(0,T]$ on $x_j(0)$, $j=1,\dots,N_1$ through bijection (\ref{eq:equivalent_odes_3}). The importance and utility of these properties for the feasibility of the proposed method as well as for the analysis of its accuracy are highlighted in the next section.

\subsection{Accuracy bounds and relevant technical assumptions}\label{sec:accuracy_bounds}

Consider system (\ref{eq:parameters_system})
\[
\begin{split}
\dot{\bfx}=&f({\bfx},\theta),\\
{\bf p}(t)=&h({\bfx}(t;\theta, {\bfx}_0),\theta), \ {\bfx}_0\in\mathbb{R}^n,
\end{split}
\]
where $f:\mathbb{R}^n\times\mathbb{R}^d\rightarrow \mathbb{R}^n$ is a piece-wise continuous and locally Lipschitz function, $h:\mathbb{R}^n\times\Real^d \rightarrow \mathbb{R}^m$ is a continuous locally Lipschitz function, $\theta\in\mathbb{R}^d$ is a vector of parameters,  $\bfx_0\in\mathbb{R}^n$ is the vector of initial conditions. Let $\Theta=({\bfx}_0,\theta)$ denote the vector of parameters and initial conditions of (\ref{eq:parameters_system}). 

Suppose that the following assumption holds true

\begin{assume}{[General properties of system (\ref{eq:parameters_system})]}\label{assume:general}
\begin{itemize}
\item $\bfx_0\in\Omega_x$, $\theta\in\Omega_\theta$,  $\Theta\in\Omega=\Omega_x\times\Omega_\theta$, where $\Omega_x$, $\Omega_\theta$, and $\Omega$ are closed and bounded subsets of $\Real^n$, $\Real^d$, and $\Real^{n+d}$, respectively 
\item  $T\in\Real$, $T>0$ is such that solutions of 
\[
\dot\bfx = f(\bfx,\theta), \ \bfx(0)=\bfx_0
\]
are defined on the interval $[0,T]$ for all $\theta\in\Omega_\theta$ and $\bfx_0\in\Omega_x$, and $\Pi \times [0,T]$ is a closed domain of $\Real^n\times\Real$ containing the corresponding solutions
\item $L_x$ and $L_\theta$ are the Lipschitz constants for the map $f$ restricted to $\Pi\times\Omega_\theta$:
\[
\begin{split}
& \|f({\bfx}, \theta) - f(\bfx',\theta)\|\leq L_x \|{\bfx} - {\bfx'}\|\\
& \|f({\bfx}, \theta) - f(\bfx,\theta')\|\leq L_\theta \|\theta - \theta'\|
\end{split}
\]
for all $(\bfx,\theta), (\bfx,\theta')$, and $(\bfx',\theta)$ in $\Pi\times\Omega_\theta$ 

\item $L_h$, $L_{h,\theta}$ are the Lipschitz constants for the maps $h(\cdot,\theta)$ and $h(\bfx,\cdot)$ restricted to $\Pi$, respectively.

\end{itemize}
\end{assume}

We note that constants $L_x$, $L_\theta$, $L_h$, $L_{h,\theta}$ always exist as long as $\Pi$ is bounded as the maps $f$ and $h$ are locally Lipschitz. In what follows, for the sake of brevity, we will denote
\[
\tilde h({\bfx}(t; \Theta),\Theta)=h({\bfx}(t; \theta,{\bfx}_0),\theta).
\]

In addition to general standard requirements for the parameter-generation system (\ref{eq:parameters_system}) imposed by Assumption \ref{assume:general}, we need to introduce further technical requirements which would enable us to quantitatively assess and analyse the problem of extrapolation of the values of $\hat{c}_k(t)$ via the parameterised map $F$ (see Eq. (\ref{eq:family_of_distributions}).  Let $[0,\tau]$, with $\tau\in(0,T]$, be an interval on which the values of $\hat{c}_k(t)$ (solutions of (\ref{eq:Smol-trunc})), $t\in[0,\tau]$, $k=1,\dots,N$ are available and known. The set
\[
\mathcal{E}(\Theta,\tau)=\{\Theta'\in\Omega \ | \tilde h({\bfx}(t; \Theta),\Theta) = \tilde h({\bfx}(t; \Theta'),\Theta') \ \mbox{for all} \ t\in[0,\tau]\}
\]
denotes the set of parameters and initial conditions of system (\ref{eq:parameters_system}) which are indistinguishable from observations $\tilde h({\bfx}(t;\Theta),\Theta)$ over the interval $[0,\tau]$. Depending on the properties of the right-hand side of (\ref{eq:parameters_system}), the set $\mathcal{E}(\Theta,\tau)$ may be a singleton, finite, or infinite. The first two cases are often referred to as {\it identifiable} in control and identification literature \cite{distefano1980parameter}, and the latter case is usually considered {\it unidentifiable}.

In our work, we do not wish to impose identifiability constraints on parameter-generation systems (\ref{eq:parameters_system}). We would, however, like to request the possibility to infer some information about the proximity of a given $\Theta$ to $\mathcal{E}(\Theta',\tau)$ from the values of $\tilde h(\bfx(t;\Theta),\Theta)$ and $\tilde h(\bfx(t;\Theta'),\Theta')$ over the interval $[0,\tau]$. This requirement is formally captured in Assumption \ref{assume:parameters} below. 

\begin{assume}{[Parameter inferability]}\label{assume:parameters} Consider (\ref{eq:parameters_system}) satisfying Assumption \ref{assume:general}. Pick $\tau\in(0,T]$. Then for any $\Theta\in\Omega$  there exists a strictly monotone function $\beta_{\theta,\tau}:\Real_{\geq 0} \rightarrow \Real_{\geq 0}$ such that for all $\Theta'\in\Omega$ the following holds true:
\[
\|\tilde h(\bfx(\cdot;\Theta),\Theta)-\tilde h(\bfx(\cdot;\Theta'),\Theta')\|_{\infty,[0,\tau]}=\max_{t\in [0,\tau]} \|\tilde h(\bfx(t;\Theta),\Theta)-\tilde h(\bfx(t;\Theta'),\Theta')\| \geq \beta_{\theta,\tau}(\|\Theta - \mathcal{E}(\Theta',\tau)\|),
\]
where 
\[
\|\Theta - \mathcal{E}(\Theta',\tau)\| = \min_{\Theta''\in\mathcal{E}(\Theta',\tau)} \|\Theta-\Theta''\|.
\]
\end{assume}

 Assumption \ref{assume:parameters} specifies the possibility to infer the distance $\|\Theta - \mathcal{E}(\Theta',\tau)\|$ from the values of the function $\tilde h(\bfx(\cdot;\Theta),\Theta)-\tilde h(\bfx(\cdot;\Theta'),\Theta')$ observed over $[0,\tau]$, $\tau\in(0,T]$ in terms of the function $\beta_{\theta,\tau}$. We note that the function may generally depend on $\Theta$ and $\tau$.
 
 The strict monotonicity of $\beta_{\theta,\tau}$ is needed for the existence of the inverse $\beta_{\theta,\tau}^{-1}$, the property we exploit in Theorem \ref{thm:exptrapolation:general}, the main statement of the section. We note that the lack of strict monotonicity of the function $\beta_{\theta,\tau}$ is not critical for producing error estimates.  If  $\beta_{\theta,\tau}$ is a mere non-decreasing function, then the inverse $\beta_{\theta,\tau}^{-1}$  can be replaced with 
 a pseudo-inverse  
 \[
 \beta^{\dagger}_{\theta,\tau}(s)=\inf\{x\in\Real_{\geq 0} | \ \beta_{\theta,\tau}(x)\geq s\}
 \]
in all relevant statements below. 

Our final set of assumptions captures the regularity constraints we impose on the parameterised family $F$ (scaling functions (\ref{eq:family_of_distributions})). The first such constraint is the requirement that the map $F(k,\cdot)$ is Lipschitz, uniformly in $k$, which is formally stated in Assumption \ref{assume:regularity:F}. 

\begin{assume}[Regularity of the map $F$]\label{assume:regularity:F} The map $F(k,\cdot)$ is Lipschitz  with the corresponding Lipschitz constant, $L_p$, uniformly in $k$. 
\end{assume}

The second regularity constraint captures the feasibility of {\it extrapolation} and {\it approximation}. To formalise this second constraint, we introduce   $\mathcal{H}_{\tau}\subset \mathcal{C}^0([0,\tau],\Real^m)$, the class of functions containing ``ideal'' parametrisations $\bfp$  of the map $F(\cdot,\bfp(t))$,  $t\in[0,\tau]$. In principle, the class $\mathcal{H}_{\tau}$ could be defined as
\begin{equation}\label{eq:parametrisation_class:h_exact}
\mathcal{H}_{\tau}=\{ \bfp\in\mathcal{C}^0([0,\tau],\Real^m) \ | \ \bfp(t)=\tilde h(\bfx(t;\Theta),\Theta) \ \mbox{for all} \ t\in[0,\tau], \ \Theta\in\Omega \},
\end{equation}
or, alternatively, as
\begin{equation}\label{eq:parametrisation_class:h_big}
\mathcal{H}_{\tau}=\{ \bfp\in\mathcal{C}^{1}([0,\tau],\Real^m) \ | \  \|\dot{\bfp}\| \leq L_h \max_{(\bfx,\theta)\in\Pi\times\Omega_\theta} \|f(\bfx,\theta)\|\}.
\end{equation}

One can now define the set $\mathcal{P} (\mathcal{H}_\tau)$ of best parametrisations $\bfp$ constrained to the class $\mathcal{H}_\tau$
\[
\mathcal{P} (\mathcal{H}_\tau)=\{ \bfp\in\mathcal{H}_\tau \ | \ \bfp=\arg \min_{v\in\mathcal{H}_\tau} \max_{k\in\Natural, k\leq N} \|F(k,v(\cdot)) - \hat{c}_k(\cdot)\|_{\infty,[0,\tau]}\}.
\]
For any $\tau\in(0,T]$, let $\Delta_\tau\in\Real_{\geq}$ be a non-negative number satisfying
\begin{equation}\label{eq:loosness_extrapolation:1}
\| \bfp' - \bfp'' \|_{\infty,[0,\tau]} \leq \Delta_\tau \ \mbox{for all} \ \bfp'\in\mathcal{P}(\mathcal{H}_\tau), \ \bfp''\in\mathcal{P}(\mathcal{H}_T).
\end{equation}
In addition, let $\gamma_\tau\in\Real_{\geq 0}$, $\tau\in(0,T]$ be such that
\begin{equation}\label{eq:loosness_extrapolation:2}
\|\tilde h(\bfx(\cdot,\Theta'),\Theta')-\tilde h(\bfx(\cdot,\Theta''),\Theta'')\|_{\infty,[0,T]} \leq \gamma_\tau \ \mbox{for all} \ \Theta'\in \Omega, \ \Theta''\in \mathcal{E}(\Theta',\tau).
\end{equation}

The values of $\Delta_\tau$, $\gamma_\tau$ characterise relevant looseness of {\it extrapolation} for functions from $\mathcal{P}(\mathcal{H}_T)$ and maps $\tilde h(\bfx(\cdot,\Theta),\Theta)$ on the basis of information contained in the behavior of these functions on $[0,\tau]$.

Similarly, for any $\tau\in(0,T]$, let $\delta_{\tau}\in\Real_{\geq 0}$ be such that
\begin{equation}\label{eq:loosness_approximation}
\max_{\bfp\in\mathcal{P}(\mathcal{H}_T)}\min_{\Theta\in\Omega} \|\bfp(\cdot) - \tilde h(\bfx(\cdot,\Theta),\Theta)\|_{\infty,[0,\tau]}\leq \delta_\tau,
\end{equation}
where $\delta_\tau$ can be viewed as a bound on the worst-case looseness of {\it approximation} of elements from $\mathcal{P}(\mathcal{H}_T)$ by  $\tilde h(\bfx(\cdot,\Theta),\Theta)$. 

Finally, let $\epsilon_{\mathrm{opt}}$ be a bound on the best admissible worst-case error between $F(k,\bfp(\cdot))$  and  $\hat{c}_k(\cdot)$: 
\begin{equation}\label{eq:error_optimal}
\epsilon_{\mathrm{opt}}= \max_{k} \|F(k,\bfp(\cdot)) - \hat{c}_k(\cdot)\|_{\infty,[0,T]}, \ \bfp\in\mathcal{P}(\mathcal{H}_T).
\end{equation}
Note that the value of $\epsilon_{\mathrm{opt}}$ in (\ref{eq:error_optimal}) is independent on the choice of $\bfp$ from the class $\mathcal{P}(\mathcal{H}_T)$.

With all the relevant assumptions and notations in place, we are now ready to state the following result.

\begin{thm}\label{thm:exptrapolation:general}[Error bounds, general] Consider system (\ref{eq:parameters_system}) and suppose that Assumptions \ref{assume:general}, \ref{assume:parameters}, \ref{assume:regularity:F} hold true. 
Let $\tau\in(0,T]$, and $\Delta_\tau$, $\gamma_\tau$, $\delta_\tau$, and $\epsilon_{\mathrm{opt}}$ be such that (\ref{eq:loosness_extrapolation:1}) -- (\ref{eq:error_optimal}) hold. Finally, let 
$\hat{\Theta}\in\Omega$ be any element of the set 
$\Omega$ satisfying
\[
\|\bfp^\ast - \tilde h(\bfx(\cdot;\hat{\Theta}),\hat{\Theta})\|_{\infty,[0,\tau]}\leq \epsilon_\tau
\]
for some $\epsilon_\tau>0$ and $\bfp^\ast\in\mathcal{P}(\mathcal{H}_\tau)$.

Then
\begin{equation}\label{eq:accuracy:general}
\begin{split}
&\|F(k,h(\bfx(\cdot;\hat\Theta),\hat{\Theta}))-\hat{c}_k(\cdot)\|_{\infty,[0,T]} \leq \\ 
& L_p \beta^{-1}_{\theta,\tau}(\epsilon_\tau + \Delta_\tau + \delta_T) \left(L_h\left(\frac{\max\{L_x,L_\theta\}}{L_x}(\exp(L_x T) -1) + 1\right)+L_{h,\theta}\right) + L_p(\gamma_\tau+\delta_T) + \epsilon_{\mathrm{opt}}
\end{split}
\end{equation}
for all $k$.
\end{thm}

{\it Proof of Theorem \ref{thm:exptrapolation:general}.} The proof of the theorem is organised into two parts. In the first part, we consider the set
\[
\Omega^\ast=\{\Theta\in\Omega \ | \  \max_{\bfp\in\mathcal{P}(\mathcal{H}_T)} \|\bfp(\cdot) - \tilde h(\bfx(\cdot,\Theta),\Theta)\|_{\infty,[0,T]}\leq \delta_T \}.
\]
This set is not empty as our assumptions require that (\ref{eq:loosness_approximation}) holds for all $\tau\in(0,T]$, including $\tau=T$. Let
\[
\mathcal{E}(\Omega^\ast,\tau)=\bigcup_{\Theta\in\Omega^\ast}\mathcal{E}(\Theta,\tau).
\]

We then produce an upper bound on 
\[
\|\hat{\Theta}- \mathcal{E}(\Omega^\ast,\tau)\|=\min_{\Theta\in \mathcal{E}(\Omega^\ast,\tau)} \|\hat{\Theta}-\Theta\|
\]
in terms of $\epsilon_\tau$, $\Delta_\tau$,  and $\delta_T$.

In the second part of the proof, we use this bound to produce the desired error estimate.

{\it Part 1}. Let $\bfp^\ast_\tau\in\mathcal{P}(\mathcal{H}_\tau)$, $\bfp^\ast_{T}\in\mathcal{P}(\mathcal{H}_T)$, and consider 
\[
\|\tilde h(\bfx(\cdot;\hat\Theta),\hat{\Theta})-\bfp_{\tau}^\ast\|_{\infty,[0,\tau]}.
\]
According to the triangle inequality
\[
\begin{split}
 \|\tilde h(\bfx(\cdot;\hat\Theta),\hat\Theta)-\bfp_{\tau}^\ast\|_{\infty,[0,\tau]}=& \|\tilde h(\bfx(\cdot;\hat\Theta),\hat\Theta)-\bfp_{T}^\ast + \bfp_{T}^\ast -\bfp_{\tau}^\ast\|_{\infty,[0,\tau]}\\
& \geq \|\tilde h(\bfx(\cdot;\hat\Theta),\hat\Theta)-\bfp_{T}^\ast\|_{\infty,[0,\tau]} - \|\bfp^\ast_{\tau}-\bfp^\ast_{T}\|_{\infty,[0,\tau]}\\
&= \|\tilde h(\bfx(\cdot;\hat\Theta),\hat\Theta)-\bfp_{T}^\ast\|_{\infty,[0,\tau]} - \Delta_\tau,
\end{split}
\]
where the last inequality follows from (\ref{eq:loosness_extrapolation:1}). 

Let $\Theta$ be any element of the set $\Omega^\ast$. Consider
\[
\begin{split}
&\|\tilde h(\bfx(\cdot;\hat\Theta),\hat\Theta)-\bfp_{T}^\ast\|_{\infty,[0,\tau]} = \|\tilde h(\bfx(\cdot;\hat\Theta),\hat\Theta) - \tilde h(\bfx(\cdot;\Theta),\Theta) + \tilde h(\bfx(\cdot;\Theta),\Theta) -\bfp_{T}^\ast\|_{\infty,[0,\tau]}\\
& \geq \|\tilde h(\bfx(\cdot;\hat\Theta),\hat\Theta) - \tilde h(\bfx(\cdot;\Theta),\Theta)\|_{\infty,[0,\tau]} - \delta_T,
\end{split}
\]
where the last inequality follows from the fact that $\Theta\in\Omega^\ast$,
\[
\|\tilde h(\bfx(\cdot;\Theta),\Theta)-\bfp_T^\ast\|_{\infty,[0,T]} \leq  \max_{\bfp\in\mathcal{P}(\mathcal{H}_T)} \|\bfp(\cdot) - \tilde h(\bfx(\cdot,\Theta),\Theta)\|_{\infty,[0,T]} \leq \delta_T,
\]
and that 
\[
\|\tilde h(\bfx(\cdot;\Theta),\Theta)-\bfp_T^\ast\|_{\infty,[0,\tau]}\leq \|\tilde h(\bfx(\cdot;\Theta),\Theta)-\bfp_T^\ast\|_{\infty,[0,T]},
\]
since $\tau \in (0,T]$.

Hence, combining the above with Assumption \ref{assume:parameters} we obtain
\[
\epsilon_\tau + \Delta_\tau + \delta_T \geq \|\tilde h(\bfx(\cdot;\hat\Theta),\hat\Theta) - \tilde h(\bfx(\cdot;\Theta),\Theta)\|_{\infty,[0,\tau]}\geq \beta_{\theta,\tau}(\|\hat\Theta-\mathcal{E}(\Theta,\tau)\|)
\]
implying that
\[
\|\hat\Theta-\mathcal{E}(\Theta,\tau)\|\leq \beta^{-1}_{\theta,\tau}(\epsilon_\tau + \Delta_\tau + \delta_T) \ \mbox{for all} \ \Theta\in\Omega^\ast,
\]
and hence
\begin{equation}\label{eq:parameters_bound}
\|\hat{\Theta}- \mathcal{E}(\Omega^\ast,\tau)\| \leq \beta^{-1}_{\theta,\tau}(\epsilon_\tau + \Delta_\tau + \delta_T).
\end{equation}

{\it Part 2}. Let 
\begin{equation}\label{eq:optimal_theta_tau}
\Theta^\ast=\arg \min_{\Theta\in\mathcal{E}(\Omega^\ast,\tau)} \|\hat{\Theta}-\Theta\|
\end{equation}
and $\theta^\ast\in\Omega_\theta$, $\bfx_0^\ast\in\Omega_x$ and $\hat{\theta}\in\Omega_\theta$, $\hat{\bfx}_0\in\Omega_x$ be such  that $\Theta^\ast=(\bfx_0^\ast,\theta^\ast)$, $\hat{\Theta}=(\hat{\bfx}_0,\hat{\theta})$.

Consider the function 
\[
\xi(t)=\bfx(t;\hat{\theta},\hat{\bfx}_0)-(\hat{\bfx}_0-\bfx_0^\ast).
\]
Given that $\bfx(\cdot;\hat{\theta},\hat{\bfx}_0)$ is defined on $[0,T]$ (see Assumption \ref{assume:general}, condition two), the function $\xi$ is defined on $[0,T]$. Moreover, as the function $f$ in (\ref{eq:parameters_system})  is locally Lipschitz, $\xi$ is the unique solution of the following initial value problem:
\[
\dot{\xi}=f(\xi + (\hat{\bfx}_0-\bfx_0^\ast),\hat{\theta}), \ \xi(0)=\bfx_0^\ast.
\]
Note that since $\hat{\bfx}_0,\bfx^\ast_0\in\Omega_x$ and $\hat{\theta},\theta^\ast\in\Omega_\theta$,  Assumption \ref{assume:general} (the second condition) implies that:
\[
\xi(t) + (\hat{\bfx}_0-\bfx_0^\ast)=\bfx(t;\hat{\theta},\hat{\bfx}_0)\in \Pi  \ \mbox{and} \ \bfx(t;\theta^\ast,\bfx_0^\ast)\in\Pi \ \mbox{for all} \ t\in [0,T].
\]

Furthermore, for any $t\in(0,T]$ we have
\[
\begin{split}
&\|\xi(t)-\bfx(t;\theta^\ast,\bfx_0^\ast)\|=\|\int_{0}^t f(\xi(s) +(\hat{\bfx}_0-\bfx_0^\ast),\hat\theta) - f(\bfx(s;\theta^\ast,\bfx_0^\ast),\theta^\ast) ds \|\\
\leq &\int_{0}^t \|f(\xi(s) +(\hat{\bfx}_0-\bfx_0^\ast),\hat\theta) - f(\bfx(s;\theta^\ast,\bfx_0^\ast),\hat\theta) + f(\bfx(s;\theta^\ast,\bfx_0^\ast),\hat\theta)  - f(\bfx(s;\theta^\ast,\bfx_0^\ast),\theta^\ast) \|ds\\
\leq & \int_0^t L_x\|\xi(s)-\bfx(s;\theta^\ast,\bfx_0^\ast)\| + L_x \|\hat{\bfx}_0-\bfx_0^\ast\| + L_\theta \| \hat{\theta}-\theta^\ast\|  ds\\
\leq & \int_0^t L_x\|\xi(s)-\bfx(s;\theta^\ast,\bfx_0^\ast)\| + \max\{L_x,L_\theta\} \|\hat\Theta-\Theta^\ast\|  ds,
\end{split}
\]
where the second inequality follows from the fact that $\xi(s) +(\hat{\bfx}_0-\bfx_0^\ast)\in\Pi$, $\bfx(s;\theta^\ast,\bfx_0^\ast)\in\Pi$ for all  $s\in [0,T]$ and that $f$ is Lipschitz in $\Pi\times\Omega_\theta$ (Assumption \ref{assume:general}, third condition).

Recall that any continuous scalar function satisfying 
\[
u(t)\leq \int_{0}^t a u(s) + b ds, \ a,b\in\Real, \ a,b>0
\]
for all $t\in[0,T]$ can be bounded as 
\[
u(t)\leq \frac{b}{a}(\exp(a t)-1).
\]
Therefore
\[
\|\xi(t)-\bfx(t;\theta^\ast,\bfx_0^\ast)\|\leq \frac{\max\{L_x,L_\theta\} \|\hat{\Theta}-\Theta^\ast\|}{L_x} (\exp(L_x t) -1)
\]
and consequently
\[
\begin{split}
\|\hat{\bfx}(t;\hat{\theta},\hat{\bfx}_0)-\bfx(t;\theta^\ast,\bfx_0^\ast)\|&\leq  \frac{\max\{L_x,L_\theta\} \|\hat{\Theta}-\Theta^\ast\|}{L_x} (\exp(L_x t) -1) + \|\hat{\bfx}_0-\bfx^\ast_0\|\\
& \leq  \frac{\max\{L_x,L_\theta\} \|\hat{\Theta}-\Theta^\ast\|}{L_x} (\exp(L_x t) -1) + \|\hat{\Theta}-\Theta^\ast\|\\
&= \|\hat{\Theta}-\Theta^\ast\| \left(\frac{\max\{L_x,L_\theta\}}{L_x}(\exp(L_x t) -1) + 1\right).
\end{split}
\]
This together with (\ref{eq:parameters_bound}) implies
\begin{equation}\label{eq:solutions_bound:parameters}
\|\hat{\bfx}(t;\hat\Theta)-\bfx(t;\Theta^\ast)\|\leq  \beta^{-1}_{\theta,\tau}(\epsilon_\tau + \Delta_\tau + \delta_T)\left(\frac{\max\{L_x,L_\theta\}}{L_x}(\exp(L_x t) -1) + 1\right) \ \mbox{for all} \ t\in[0,T].
\end{equation}
Finally, consider
\[
\begin{split}
& \|F(k,\tilde h(\bfx(\cdot;\hat\Theta),\hat{\Theta}))-\hat{c}_k(\cdot)\|_{\infty,[0,T]}=\|F(k,\tilde h(\bfx(\cdot;\hat\Theta),\hat{\Theta}))-F(k,\bfp^\ast_{T}(\cdot)) + F(k,\bfp^\ast_{T}(\cdot))-\hat{c}_k(\cdot)\|_{\infty,[0,T]}\\
& \leq \|F(k,\tilde h(\bfx(\cdot;\hat\Theta),\hat{\Theta}))-F(k,\bfp^\ast_{T}(\cdot))\|_{\infty,[0,T]} + \epsilon_{\mathrm{opt}}\leq L_p\|\tilde h(\bfx(\cdot;\hat\Theta),\hat{\Theta}) - \bfp^\ast_{T}(\cdot)\|_{\infty,[0,T]}+\epsilon_{\mathrm{opt}},
\end{split}
\]
where the last inequality is due to  (\ref{eq:accuracy:general}) and the fact that the map $F(k,\cdot)$ being Lipschitz uniformly in $k$. 

Taking into account that $\hat{\bfx}(t;\hat{\theta},\hat\bfx_0), \bfx(t;\theta^\ast,\bfx_0^\ast)$ are in $\Pi$ for all $t\in[0,T]$, and that $h$ is Lipschitz in $\Pi$ (see Assumption \ref{assume:general}), we obtain:
\begin{equation}\label{eq:pre-final-bound:1}
\begin{split}
& \|\tilde h(\bfx(\cdot;\hat\Theta),\hat{\Theta}) - \bfp^\ast_{T}(\cdot)\|_{\infty,[0,T]}\leq \|\tilde h(\bfx(\cdot;\hat\Theta),\hat{\Theta}) - \tilde h(\bfx(\cdot;{\Theta}^\ast),{\Theta}^\ast)\|_{\infty,[0,T]} +  \|h(\bfx(\cdot;{\Theta}^\ast),{\Theta}^\ast)- \bfp^\ast_{T}(\cdot)\|_{\infty,[0,T]}\\
&\leq\|\tilde h(\bfx(\cdot;\hat\Theta),\hat{\Theta}) - \tilde h(\bfx(\cdot;{\Theta}^\ast),{\Theta}^\ast)\|_{\infty,[0,T]} +  \|\tilde h(\bfx(\cdot;{\Theta}^\ast),{\Theta}^\ast)- \tilde h(\bfx(\cdot;{\Theta}),{\Theta})\|_{\infty,[0,T]} + \delta_T,
\end{split}
\end{equation}
where $\Theta\in\Omega^\ast$ is chosen so that $\Theta^\ast\in\mathcal{E}(\Theta,\tau)$. Note that the existence of such $\Theta$ is ensured by the way $\Theta^\ast$ is defined in (\ref{eq:optimal_theta_tau}).  Observe that
\begin{equation}\label{eq:pre-final-bound:2}
\begin{split}
& \|\tilde h(\bfx(\cdot;\hat\Theta),\hat{\Theta}) - \tilde h(\bfx(\cdot;{\Theta}^\ast),{\Theta}^\ast)\|_{\infty,[0,T]}  \\
& \ \ \ \ \ \ \ \ \ \  \ \ \ \ \ \ \ \ \ \ \leq \|\tilde h(\bfx(\cdot;\hat\Theta),\hat{\Theta}) - \tilde h(\bfx(\cdot;{\Theta}^\ast),\hat{\Theta})\|_{\infty,[0,T]} + \|\tilde h(\bfx(\cdot;\Theta^\ast),\hat{\Theta}) - \tilde h(\bfx(\cdot;{\Theta}^\ast),{\Theta}^\ast)\|_{\infty,[0,T]}\\
&\ \ \ \ \ \ \ \ \ \ \ \ \ \ \ \ \ \ \ \ \leq L_h \|\bfx(\cdot;\hat\Theta)-\bfx(\cdot;\Theta^\ast)\| + L_{h,\theta}\|\hat\Theta-\hat\Theta^\ast\|.
\end{split}
\end{equation}
Moreover, recalling that $\Theta^\ast\in\mathcal{E}(\Theta,\tau)$ and taking (\ref{eq:loosness_extrapolation:2}) into account, one can conclude that 
\begin{equation}\label{eq:pre-final-bound:3}
 \|\tilde h(\bfx(\cdot;{\Theta}^\ast),{\Theta}^\ast)- \tilde h(\bfx(\cdot;{\Theta}),{\Theta})\|_{\infty,[0,T]}\leq \gamma_\tau.
\end{equation}

Given that the map $F(k,\cdot)$ is Lipschitz (Assumption \ref{assume:regularity:F}) and invoking (\ref{eq:pre-final-bound:1})--(\ref{eq:pre-final-bound:3}) we arrive at the following estimate: 
\[
\|F(k,\tilde h(\bfx(\cdot;\hat\Theta),\hat{\Theta}))-\hat{c}_k(\cdot)\|_{\infty,[0,T]}\leq L_p L_h \|\bfx(\cdot;\hat\Theta)-\bfx(\cdot;\Theta^\ast)\| + L_p L_{h,\theta}\|\hat\Theta-\hat\Theta^\ast\| + L_p(\gamma_\tau+\delta_T) + \epsilon_{\mathrm{opt}}.
\]

Combining this estimate with (\ref{eq:solutions_bound:parameters}), (\ref{eq:parameters_system}) results in
\[
\begin{split}
&\|F(k,\tilde h(\bfx(\cdot;\hat\Theta),\hat{\Theta}))-\hat{c}_k(\cdot)\|_{\infty,[0,T]} \leq \\
& L_p \beta^{-1}_{\theta,\tau}(\epsilon_\tau + \Delta_\tau + \delta_T) \left(L_h\left(\frac{\max\{L_x,L_\theta\}}{L_x}(\exp(L_x T) -1) + 1\right)+L_{h,\theta}\right) + L_p(\gamma_\tau+\delta_T) + \epsilon_{\mathrm{opt}}.
\end{split}
\]
 $\square$

\begin{rem}[Parametrisation classes $\mathcal{H}_\tau$ and knowledge-informed prediction]\label{rem:simple_bound:general} \normalfont The statement of Theorem \ref{thm:exptrapolation:general} and error bound (\ref{eq:accuracy:general}) can be simplified if $\mathcal{H}_\tau$ is defined as in (\ref{eq:parametrisation_class:h_exact}). In this case, $\delta_T=0$, and  (\ref{eq:accuracy:general}) becomes:
\begin{equation}\label{eq:accuracy:exact}
\begin{split}
&\|F(k,h(\bfx(\cdot;\hat\Theta),\hat{\Theta}))-\hat{c}_k(\cdot)\|_{\infty,[0,T]} \leq \\ 
& L_p \beta^{-1}_{\theta,\tau}(\epsilon_\tau + \Delta_\tau) \left(L_h\left(\frac{\max\{L_x,L_\theta\}}{L_x}(\exp(L_x T) -1) + 1\right)+L_{h,\theta}\right) + L_p(\gamma_\tau) + \epsilon_{\mathrm{opt}}.
\end{split}
\end{equation}
In addition, if system (\ref{eq:parameters_system}) is such that the maps $\tilde{h}(\bfx(t;\Theta'),\Theta')$ and $\tilde{h}(\bfx(t;\Theta''),\Theta'')$ coincide on $[0,T]$ whenever $\Theta'\in\Omega$ and $\Theta''\in\mathcal{E}(\Theta',\tau)$, $\tau\in(0,T]$ (that is if  $\Theta''\in\mathcal{E}(\Theta',\tau)$ implies that $\Theta''\in\mathcal{E}(\Theta',T)$) then $\gamma_\tau$ can be set to zero leading to a further simplification of  (\ref{eq:accuracy:general}), (\ref{eq:accuracy:exact}):
\[ 
\begin{split}
&\|F(k,h(\bfx(\cdot;\hat\Theta),\hat{\Theta}))-\hat{c}_k(\cdot)\|_{\infty,[0,T]} \leq \\ 
& L_p \beta^{-1}_{\theta,\tau}(\epsilon_\tau + \Delta_\tau) \left(L_h\left(\frac{\max\{L_x,L_\theta\}}{L_x}(\exp(L_x T) -1) + 1\right)+L_{h,\theta}\right)  + \epsilon_{\mathrm{opt}}.
\end{split}
\]

In the context of further simplification, we note that the inclusion of additional information about the desired behavior of functions from $\mathcal{H}_\tau$, such as the requirements that they satisfy conditions similar to (\ref{eq:params_modified_property}), has the potential to reduce the values of $\Delta_\tau$ and $\gamma_\tau$ in bounds stemming from Theorem \ref{thm:exptrapolation:general}. This, in effect, captures what we called here {\it knowledge-informed learning}: the increased accuracy of prediction guided by controlling the size of hypothesis classes $\mathcal{H}_\tau$ based on prior phenomenological knowledge about the physics of the processes being predicted.

\end{rem}

Theorem \ref{thm:exptrapolation:general} and Remark \ref{rem:simple_bound:general}, taken together with Theorem \ref{thm:neural_network:dynamic}, enable to produce an extrapolation error bound for mappings modelled by (\ref{eq:reference_net}).

\begin{corollary}[Error bounds for NNs]\label{cor:neural_net_error} Consider the class of neural networks defined by (\ref{eq:reference_net}) with 
\[
\begin{split}
&W^2\in [-a_1,a_1]^{N_1}, \ W^1\in [-a_2,a_2]^{N_1}, \ b^2\in[-a_3, a_3], \ b^1\in[-a_4,a_4]^{N_1}, \ a_1,a_2,a_3,a_4\in \Real_{>0},\\
&\Omega_{W,b}=[-a_2,a_2]^{N_1}\times[-a_1,a_1]^{N_1}\times[-a_3,a_3]\times[-a_4,a_4]^{N_1},
\end{split}
\]
and the corresponding system of ordinary differential equations (\ref{eq:equivalent_odes_2}) with $\bfx=(x_1,\dots,x_{N_1})$, $\theta=(W^2,W^1,b^2,b^1)$, initial conditions 
\[
\bfx_0=(x_1(0),\dots,x_{N_1}(0)), \ x_j(0)=\frac{1}{1+\exp(-b^1_j)}, \ j=1,\dots,N_1,
\]
appropriately chosen $\Omega_\theta$, $\Omega_x$, and mapping $h$:
\[
h(\bfx(t;\theta,\bfx_0(b^1)),\theta)=(W^2,\bfx(t;\theta,\bfx_0(b^1))+b^2.
\]
Let $\tau\in(0,T]$, and suppose that Assumptions \ref{assume:parameters}, \ref{assume:regularity:F} hold true. Furthermore, let $\mathcal{H}_\tau$ be defined by (\ref{eq:parametrisation_class:h_exact}), and $\Delta_\tau$, $\gamma_\tau$, and $\epsilon_{\mathrm{opt}}$ be defined as in (\ref{eq:loosness_extrapolation:1}), (\ref{eq:loosness_extrapolation:2}), and (\ref{eq:error_optimal}).  

Let ${W^{\ast}}^{2}$, ${W^{\ast}}^{1}$, ${b^{\ast}}^{2}$, ${b^{\ast}}^1$ be parameters of (\ref{eq:reference_net}) minimising
\[
\|F(k,y(\cdot, W^2, W^1, b^2, b^1)) - \hat{c}_k(\cdot)\|_{\infty,[0,\tau]}
\]
over $\Omega_{W,b}$, and $(\hat{W}^{2}, \hat{W}^{1}, \hat{b}^{2}, \hat{b}^{1})\in\Omega_{W,b}$ be such that
\[
\|y(\cdot, {W^{\ast}}^2, {W^{\ast}}^{1}, {b^{\ast}}^2, {b^{\ast}}^{1}) - y(\cdot, \hat{W}^2, \hat{W}^1, \hat{b}^2, \hat{b}^1)\|_{\infty,[0,\tau]}\leq \epsilon_\tau.
\]

Then
\begin{equation}\label{eq:accuracy:neural_net}
\begin{split}
&\|F(k,y(\cdot, \hat{W}^2, \hat{W}^1, \hat{b}^2, \hat{b}^1))-\hat{c}_k(\cdot)\|_{\infty,[0,T]} \leq \\ 
& L_p \beta^{-1}_{\theta,\tau}(\epsilon_\tau + \Delta_\tau) \left(a_1\sqrt{N_1} \left(\frac{\max\{a_2,1/4\}}{a_2}(\exp(a_2 \sqrt{N_1} T) -1) + 1\right)+\sqrt{N_1+1}\right) + L_p(\gamma_\tau) + \epsilon_{\mathrm{opt}}.
\end{split}
\end{equation}
\end{corollary}
{\it Proof of Corollary \ref{cor:neural_net_error}}. The corollary follows immediately from Theorem \ref{thm:exptrapolation:general} and Remark \ref{rem:simple_bound:general}. Indeed, observe that solutions of (\ref{eq:equivalent_odes_2}) are defined for all $t\in\Real$, are bounded and remain in $[0,1]^{N_1}$ 
for all $W^1\in\Real^{N_1}$, $x_j(0)=1/(1+\exp(-b_j^1))$, $b_j^1\in[-a_4,a_4]$. 

Define the set
\[
\Pi=\{\bfx\in\Real^{N_1} \ | x_j\in[0,1], \ j=1,\dots,N_1\}.
\]

The right-hand side of system (\ref{eq:equivalent_odes_2}) is Lipschitz in $\Pi\times\Omega_{\theta}$ with
\[
L_x= a_2 \sqrt{N_1}, \ L_\theta= \frac{1}{4}\sqrt{N_1}, 
\]
where the factor $1/4$ in $L_\theta$ is due to that $\max_{x\in[0,1]} x(1-x)=1/4$. The map $h$ is also Lipschitz with
\[
L_h=a_1\sqrt{N_1}, \ L_{h,\theta}=\sqrt{N_1+1}.
\]
Therefore Assumption \ref{assume:general} holds true, and  bound (\ref{eq:accuracy:neural_net}) follows from (\ref{eq:accuracy:exact}). $\square$

\begin{rem}[Networks with other activation functions] \normalfont Results and the approach presented so far may be extended to networks with other activation functions $\sigma$, such as the hyperbolic tangent $\tanh$, as long as they admit an equivalent representation by a system of ordinary differential equations. 
\end{rem}

\subsection{Training protocol and meta-parameters}\label{sec:training_protocols}

In the paper we investigated the application of the proposed approach to the aggregation process, which was defined by the system:

\begin{align}
    & \frac{d \hat{c}_k}{dt} = \frac{1}{2} \sum_{i+j=k} K_{ij} \hat{c}_i \hat{c}_j - \hat{c}_k \sum_{j=1}^{N} K_{kj} \hat{c}_j + \delta_{k,1},\qquad k = 1..N, \nonumber \\ 
    & \hat{c}_k(0) = \delta_{k,1} \nonumber 
\end{align}

where aggregation kernel $K_{ij}$ took three different forms:

\begin{align}
    & K_{ij} = 1, \ \text{ unit kernel} \nonumber \\
    & K_{ij} = (i \cdot j)^{0.2}, \ \text{ product kernel} \nonumber \\
    & K_{ij} = i^{0.5} + j^{0.5}, \ \text{ sum-like kernel} \nonumber \\
\end{align}

For each kernel aggregation process was modelled with fast-solver based on low-rank approximations described in \cite{matveev2014fasttransliteration, matveev2015fast}. Each system had 40000 equations, that is $N = 40000$. Each process was calculated on a time-grid with $N_T = 8000$ temporal steps. Duration of a single step $dt$ differed from kernel to kernel, because processes had different characteristic time. Hence, the data of an aggregation process was a matrix $A \in \Real^{40000 \times 8000}$.

All the calculations and training presented below were done in Matlab (version R2022a) on Intel(R) Core(TM) i7-7700HQ.

The calculated solution was modified with the logarithmic transformation with cut-off \ref{eq:solution-transformation} to avoid errors related to machine precision.
The first step of knowledge-informed prediction is retrieving parameters $\boldsymbol{p}(t) = (W(t), B(t))$ from the parametrizing network $NN_p$. The network had the following structure \\
\begin{center}
\begin{tabular}{P{0.1\textwidth}|P{0.3\textwidth}|P{0.1\textwidth}}
        Layer & Type & Size\\
        \hline \hline
        1 & Input & 1 \\
        2 & Fully Connected & 1 \\
        3 & ReLU & \\
        4 & Fully Connected & 1 
\end{tabular}\\
\end{center}
where the weight parameter from the 2nd layer is frozen to be equal $-1$, and the bias parameter from the 4th layer is frozen to be equal to cut-off value of $\log(10^{-7})$.
Parameters were calculated on a rather sparse time-grid with length of the time step $dt_p = 20dt$, and then were interpolated with splines. The result of training of parametrizing network is two retrieved parameters: $W(t)$, the weight from the fourth layer, and $B(t)$, the bias from the second layer. After retrieving we applied bijective transformations $\mathcal{T}_W(\cdot)$ and $\mathcal{T}_B(\cdot)$ to the parameters.

The second step is the training of the predicting network $NN_e$ for all retrieved parameters separately. The values of parameters calculated for 8000 temporal steps were split into three sets: training (steps from 1000 to 1900), validating (steps from 1901 to 2000), and testing (steps from 2001 to 8000). Steps from 0 to 999 were omitted in the training sets, because behavior of parameters is non-monotonous on this interval. The predicting network had the following structure
\begin{center}
\begin{tabular}{P{0.1\textwidth}|P{0.3\textwidth}|P{0.1\textwidth}}
        Layer & Type & Size\\
        \hline \hline
        1 & Input & 1 \\
        2 & Fully Connected & 5 \\
        3 & Sigmoid & \\
        4 & Fully Connected & 1 
\end{tabular}\\
\end{center}
without batch normalisation and dropout layers. The two copies of the network were trained independently: one for $\mathcal{T}_W(W(t))$, transformed $W(t)$, and another for $\mathcal{T}_B(B(t))$, transformed $B(t)$. The training of the networks was aimed not only to convergence with training data, but also to ensure that networks' output obeys properties that were drawn for extrapolating functions. Hence, training loop was two-step. On the first step, networks input was an interval of $t \in [0, \tau]$, which was represented by temporal steps from 1 to 1900. The training step was made with ADAM optimizer \cite{kingma2014adam} with the loss function $L_u(NN_e,u) = \int_0^{\tau} |NN_e(t) - u(t)|^2 d t$, where $u(t)$ is a target function, which is $\mathcal{T}_B(B(t)$ or $\mathcal{T}_B(B(t)$. On the second step, the input was the whole target interval $t \in [0, T]$, which is represented by temporal steps from 1 to 8000. Here we also used ADAM optimizer, but with another loss function $L_d(NN_e) = \int_0^T C_1 \cdot ReLU(NN_e'(t)) + C_2 \cdot ReLU(NN_e''(t)) + C_3 \cdot ReLU(NN_e^{(3)}(t)) dt$. Note that on the second step there was no target function. The optimizer was used with the default training options: initial learning rate was equal $0.001$, gradient decay factor was equal $0.9$, squared gradient decay factor was equal $0.999$.

The training was validated with the loss function $L_u$ on temporal steps from 1901 to 2000. As the learning process occurred non-monotonically with periodical oscillations of values of the loss function $L_u$ on the training set $Loss_t$, in which the minimum of values of the loss function on the validation set $Loss_v$, we used the following stopping criteria: we exit the training loop, when $Loss_v$ reached a value less than $10^{-6}$ on an epoch in which the global minimum of $Loss_t$ was reached.

\end{document}